\documentclass[11pt]{article}

\usepackage[utf8]{inputenc}
\usepackage[T1]{fontenc}
\usepackage{lmodern}

\usepackage{amsmath,amssymb,amsfonts,amsthm,mathtools}
\usepackage{bm}
\usepackage{mathrsfs}
\usepackage{enumitem}
\usepackage{booktabs}
\usepackage{array}
\usepackage{graphicx}
\usepackage{xcolor}
\usepackage{microtype}
\usepackage{geometry}
\usepackage{hyperref}
\usepackage{tikz}
\usepackage{tikz-cd}
\hypersetup{
    colorlinks=true,
    linkcolor=blue,
    citecolor=blue,
    urlcolor=blue
}

\newtheorem{theorem}{Theorem}[section]
\newtheorem{proposition}[theorem]{Proposition}
\newtheorem{lemma}[theorem]{Lemma}
\newtheorem{corollary}[theorem]{Corollary}

\theoremstyle{definition}
\newtheorem{definition}[theorem]{Definition}
\newtheorem{example}[theorem]{Example}

\theoremstyle{remark}
\newtheorem{remark}[theorem]{Remark}

\newcommand{\R}{\mathbb{R}}

\newcommand{\dd}{\mathrm{d}}
\newcommand{\Lie}{\mathcal{L}}
\newcommand{\id}{\mathrm{id}}
\newcommand{\cO}{\mathcal{O}}
\newcommand{\Z}{\mathbb{Z}}

\title{
\textbf{Structure-Preserving Integrators for Locally Conformally
Symplectic Systems: Variational Methods and Globalization}
}

\author{
Ayten Gezici$^{1,2}$\thanks{Corresponding author.
Email: \href{mailto:agezici@ucsd.edu}{\texttt{agezici@ucsd.edu} }
}
\and
Melvin Leok$^{2}$\thanks{
Email: \href{mailto:mleok@ucsd.edu}{\texttt{mleok@ucsd.edu}}
}
\and
Cristina Sard\'on$^{3}$\thanks{
Email: \href{mailto:mariacristina.sardon@upm.es}{\texttt{mariacristina.sardon@upm.es}}
}
}

\date{}

\begin{document}

\maketitle

\begin{center}
\small
$^{1}$Department of Mathematics,
Gebze Technical University,
41400 Gebze-Kocaeli, Türkiye
\\[1mm]
$^{2}$Department of Mathematics,
University of California, San Diego,
9500 Gilman Drive, La Jolla, CA 92093-0112, USA
\\[1mm]
$^{3}$Department of Applied Mathematics,
Universidad Polit\'ecnica de Madrid,
Av. Juan de Herrera 6,
28040 Madrid, Spain
\end{center}

\vspace{0.5cm}


\begin{abstract}
Locally conformally symplectic (LCS) manifolds are characterized by a
nondegenerate two-form $\Omega$ satisfying $\dd\Omega=\varphi\wedge\Omega$,
where $\varphi$ is a closed one-form, called the Lee form. Although an
LCS structure is locally conformal to an ordinary symplectic structure,
the local symplectic descriptions need not arise from a single global
symplectic form. This creates a natural local-to-global problem for
geometric numerical integration: local symplectic discretizations must
be compatible on chart overlaps in order to determine a globally
defined numerical dynamics.

We develop a Hamiltonian variational framework for the construction and
globalization of structure-preserving integrators for LCS systems. The
continuous theory is first recalled in a form adapted to a Type~II
Hamiltonian variational principle. We then introduce exact and
approximate discrete right Hamiltonians for the local symplectic
representatives and identify the transformation law that they must
satisfy on overlaps. This compatibility condition implies equivariance
of the associated local discrete Hamiltonian maps and therefore provides
a sufficient condition for their globalization. The resulting global
discrete map preserves the LCS two-form conformally, with a discrete
conformal factor determined by the local conformal potentials.

The framework applies naturally to Type~II Hamiltonian variational
integrators and, in the cotangent-type setting considered here, to
symplectic Runge--Kutta methods. In particular, Gauss--Legendre methods
and the implicit midpoint rule yield globally compatible LCS
structure-preserving schemes. Under suitable regularity assumptions,
the corresponding Lagrangian formulation is recovered through discrete
Legendre transformations. Numerical examples illustrate overlap
compatibility, chart switching, convergence, and discrete conformal
preservation.
\end{abstract}

\medskip



\section{Introduction}
\label{sec:introduction}

Geometric numerical integration is concerned with the construction of
numerical methods that reproduce, at the discrete level, geometric
properties of the differential equations being approximated. In
Hamiltonian mechanics, the central example is provided by symplectic
integrators, which preserve the symplectic structure of phase space and
typically exhibit favorable long-time behavior that is not shared by
generic numerical schemes. We refer to \cite{Blanes-Casas,HaLuWa, Reich-96, Sanz-Serna} for more details about structure-preserving integration techniques. Variational integrators provide a closely
related construction based on a discretization of Hamilton's principle \cite{DuruisseauxLeok22, Farr-2009, MarsWest},
while Hamiltonian variational integrators arise from discrete
Hamiltonian variational principles and generating functions \cite{DuSchmittLeok,LallWest-discrete,LeZh11}.

Locally conformally symplectic geometry provides a natural extension of
the symplectic setting. An LCS manifold is a triple
$(M,\Omega,\varphi)$ consisting of a smooth manifold $M$, a
nondegenerate two-form $\Omega$, and a closed one-form $\varphi$,
called the Lee form, satisfying
\begin{equation}
    \dd\Omega
    =
    \varphi\wedge\Omega.
    \label{eq:intro-LCS}
\end{equation}
Since $\varphi$ is closed, it is locally exact,
$\varphi|_{U_\alpha}=\dd\sigma_\alpha$ on a suitable open covering
$\{U_\alpha\}_{\alpha\in A}$ of $M$, and
$\Omega_\alpha=e^{-\sigma_\alpha}\Omega|_{U_\alpha}$ is closed and
nondegenerate, i.e., an ordinary symplectic form on $U_\alpha$. For
more details on LCS manifolds, see
\cite{Banyaga,Bazzoni2018,LeeHC,vaisman}. Hamilton--Jacobi theory in this setting has been developed in \cite{esen21lcs}.

This local symplectic description suggests an apparently simple
strategy for numerical integration: construct a symplectic integrator
in each local symplectic representation and then reconstruct a
numerical method for the original LCS system. The essential difficulty,
however, is global. A collection of locally defined symplectic
numerical maps does not, by itself, determine a global discrete
dynamics. The local maps must agree on overlaps after the corresponding
changes of representation are taken into account. Indeed, on a
connected component of an overlap $U_\alpha\cap U_\beta$ the potentials
differ by a constant, $\sigma_\alpha-\sigma_\beta=C_{\alpha\beta}$, so
that $\Omega_\beta=c_{\alpha\beta}\Omega_\alpha$ with
$c_{\alpha\beta}=e^{C_{\alpha\beta}}$, and the numerical constructions
performed in the two local symplectic representations must be
compatible with this overlap structure.

The purpose of the present work is to formulate this compatibility
problem at the level of discrete Hamiltonian mechanics and to provide
conditions under which local structure-preserving discretizations
globalize to a well-defined LCS numerical dynamics. Our starting point
is Hamiltonian rather than Lagrangian. This choice is important because
the Hamiltonian formulation does not require the hyperregularity
assumptions needed to pass between Hamiltonian and Lagrangian
descriptions. We first recall the continuous LCS Hamiltonian theory and
formulate its local canonical dynamics through a Type~II Hamiltonian
variational principle \cite{LeZh11,ScLe18}; this
continuous material fixes the geometric and variational notation used
in the discrete construction.

The discrete theory begins with the exact discrete right Hamiltonian,
the Type~II generating function of the exact time-$h$ flow of a local
Hamiltonian representative $H_\alpha$. Replacing it by a computable
approximation produces a local discrete Hamiltonian map which, as in
ordinary discrete Hamiltonian mechanics, is symplectic in the local
canonical variables. Local symplecticity, however, does not solve the
globalization problem: the question is whether the discrete
constructions obtained in two overlapping local representations
describe the same numerical trajectory. In the cotangent-type local representation considered in this paper, we first consider overlap transformations of the form
\begin{equation}
    T_{\alpha\beta}(q,r)
    =
    (q,c_{\alpha\beta}r).
    \label{eq:intro-transition}
\end{equation}
We show that the exact discrete right Hamiltonians satisfy a definite
transformation law under this map. This motivates the corresponding
compatibility condition for approximate discrete Hamiltonians,
\begin{equation}
    H_{\beta,h}^{d+}(q,s)
    =
    c_{\alpha\beta}
    H_{\alpha,h}^{d+}
    \left(
        q,c_{\alpha\beta}^{-1}s
    \right),
    \label{eq:intro-discrete-compatibility}
\end{equation}
which is a condition on the discrete generating functions themselves,
rather than an a priori assumption on the resulting numerical maps.
Under \eqref{eq:intro-discrete-compatibility} the local discrete maps
satisfy the equivariance relation
$F_{\beta,h}\circ T_{\alpha\beta}=T_{\alpha\beta}\circ F_{\alpha,h}$,
which is precisely the commutativity condition required for the local
maps to glue, and after reconstruction in the original LCS variables
the resulting global map $F_h$ satisfies the discrete conformal
preservation law
\begin{equation}
    F_h^*\Omega
    =
    e^{\sigma_\alpha\circ F_h-\sigma_\alpha}\Omega,
    \label{eq:intro-discrete-LCS}
\end{equation}
which reproduces the characteristic conformal preservation property of
the continuous LCS Hamiltonian flow. Symplecticity is thus a local
geometric property of the numerical method, whereas equivariance is the
compatibility property that allows those local methods to represent a
single global dynamics; the globalization statement itself is
straightforward once equivariance is known, and the substantive issue
is to identify verifiable conditions on the discretization that
guarantee it.

The Type~II framework provides such conditions naturally \cite{LeZh11}:
Galerkin approximations of the Hamiltonian action can be constructed so
as to inherit the required scaling law. Moreover, because the
transition maps \eqref{eq:intro-transition} are linear, standard
Runge--Kutta discretizations are equivariant when the same tableau is
used in every local representation, and if the method is additionally
symplectic the resulting global method is an LCS structure-preserving
integrator. This includes the Gauss--Legendre family and, as its
simplest member, the implicit midpoint rule \cite{HaLuWa,MarsWest}.
The restriction to the cotangent-type setting is relevant here: for
general nonlinear coordinate transformations, Runge--Kutta methods are
not automatically equivariant, and we therefore formulate the principal
results under explicit hypotheses on the local symplectic
representations and their transition maps rather than claiming
coordinate invariance for arbitrary LCS Darboux charts. Once the
Hamiltonian theory has been established, a Lagrangian counterpart can
be obtained under appropriate regularity assumptions, see
\cite{AbMa78,Marsden1999,MarsWest}; in this sense the Lagrangian
formulation is complementary to, rather than logically prior to, the
Hamiltonian globalization result.

The framework also extends to affine cotangent transitions
\[
    T_{\alpha\beta}(q,r)
    =
    \bigl(Bq+b,c_{\alpha\beta}B^{-T}r\bigr),
\]
where $B$ is constant and invertible, $b$ is constant, and
$c_{\alpha\beta}>0$ is constant on each connected overlap
component. We derive the associated Type~II transformation
law, including the endpoint term arising from the translation,
and verify the equivariance of Runge--Kutta methods
under these transitions.

We also clarify the distinction between the LCS Hamiltonian
framework and the conformal Hamiltonian systems studied in
\cite{mclachlan2001conformal}, for which conformally symplectic
integrators have been developed, for example, in
\cite{BhFlMo-2016, FrSuRo-2020}.
For Hamiltonian systems with constant damping rate $\gamma$,
the exact flow satisfies
$\psi_t^*\omega = e^{-\gamma t}\omega$,
where $\omega$ is a closed symplectic form, and conformally
symplectic integrators reproduce the corresponding discrete
scaling law. In contrast, the flow of an LCS Hamiltonian system
satisfies
\[
    \psi_t^*\Omega=\rho_t\,\Omega,
    \qquad
    \rho_t(z)=
    \exp\!\left(
        \int_0^t \varphi(\xi_H)(\psi_s(z))\,\dd s
    \right),
\]
so that the conformal factor generally depends on the
initial point, rather than being determined solely by a
constant damping rate. Even when the Lee form is exact,
$\varphi=\dd\sigma$, one obtains
$\rho_t= e^{\sigma\circ \psi_t-\sigma}$,
which does not in general equal $e^{-\gamma t}$;
in this exact case, the flow preserves the globally defined
symplectic form $e^{-\sigma}\Omega$. Our discrete construction
preserves the corresponding geometric relation
$F_h^*\Omega=e^{\sigma_\alpha\circ F_h-\sigma_\alpha}\Omega$
for steps represented within a single chart.
When the Lee form is non-exact, the conformal potentials are
available only locally, and compatibility on chart overlaps
is required. Our contribution is to construct compatible local
symplectic integrators that define a global LCS
structure-preserving method.

The results presented here build on the discrete locally conformal framework developed in \cite{EsGeGu24}. While that work derives the discrete LCS Hamiltonian equations, we recover these equations by applying a Type II Hamiltonian variational principle in each local chart and expressing the resulting discrete Hamilton equations in LCS coordinates. Our contributions include this direct Hamiltonian variational derivation and the construction of globally defined integrators through compatibility conditions for local symplectic integrators on chart overlaps. We formulate and verify these conditions and provide numerical examples illustrating convergence, preservation of the LCS conformal scaling law, and consistency under chart switching.

The paper is organized as follows.
Section~\ref{sec:continuous} recalls the continuous LCS Hamiltonian
geometry and develops the Type~II variational formulation required in
the sequel. Section~\ref{sec:discrete-ham} introduces exact and
approximate discrete right Hamiltonians and the associated local
symplectic maps. Section~\ref{sec:globalization} derives the overlap
transformation law of the exact discrete right Hamiltonian, introduces
the compatibility condition for approximate discrete right
Hamiltonians, shows that it implies equivariance of the local discrete
maps, and establishes the globalization and discrete LCS preservation
results. Section~\ref{sec:compatible-methods} identifies concrete
classes of compatible methods, including Type~II variational
constructions and symplectic Runge--Kutta methods, and extends the
relevant results to affine base transitions.
Section~\ref{sec:lagrangian} recovers the corresponding discrete
Lagrangian formulation under regularity assumptions.
Section~\ref{sec:numerics} investigates the resulting methods
numerically through convergence, overlap-consistency, chart-switching,
and structure-preservation tests, and
Section~\ref{sec:conclusions} summarizes the results and discusses an
open problem.


\section{Continuous LCS Hamiltonian mechanics and variational formulation}
\label{sec:continuous}

In this section, we recall the continuous geometric framework underlying
the discrete constructions developed later in the paper. Locally
conformally symplectic geometry and the corresponding Hamiltonian and
Lagrangian dynamics have been extensively studied in the literature
\cite{Banyaga,Bazzoni2018,esen21lcs,LeeHC,Libermann-LCS,vaisman}; the purpose of
this section is not to develop a new continuous theory, but to collect
the ingredients that will be needed for the Hamiltonian discretization
and globalization problem. The Hamiltonian viewpoint will be primary
throughout the paper, and no hyperregularity assumption is required at
this stage; the Lagrangian description is developed, under appropriate
regularity assumptions, in Section~\ref{sec:lagrangian}.

\paragraph{Locally conformally symplectic manifolds.}

\begin{definition}
A \emph{locally conformally symplectic manifold} is a triple
$(M,\Omega,\varphi)$, where $M$ is a smooth even-dimensional manifold,
$\Omega\in\Omega^2(M)$ is a nondegenerate two-form, and
$\varphi\in\Omega^1(M)$ is a closed one-form such that
\begin{equation}
    \dd\Omega
    =
    \varphi\wedge\Omega.
    \label{eq:LCS-definition}
\end{equation}
The one-form $\varphi$ is called the \emph{Lee form}.\footnote{When
$\dim M=2$, the identity $\dd\Omega=\varphi\wedge\Omega$ holds
automatically, since both sides are $3$-forms and hence vanish.
Thus, $\Omega$ alone does not determine the Lee form; a closed
one-form $\varphi$ is specified as part of the LCS triple
$(M,\Omega,\varphi)$. When $\dim M=2n\geq4$, the map
$\alpha\mapsto\alpha\wedge\Omega$ on $1$-forms is injective,
so the Lee form, when it exists, is uniquely determined
by $\Omega$.}
\end{definition}

\noindent
For the closed one-form $\varphi$, the associated
Lichnerowicz--de~Rham differential (see \cite{GuLi84}) is defined by
\begin{equation}
    \dd_\varphi
    :=
    \dd-\varphi\wedge.
    \label{eq:Lichnerowicz-differential}
\end{equation}
Since $\dd\varphi=0$, one has $\dd_\varphi^2=0$, and the LCS condition
\eqref{eq:LCS-definition} is equivalently written as
$\dd_\varphi\Omega=0$.

The terminology ``locally conformally symplectic'' follows from the
local exactness of the Lee form. Let $\{U_\alpha\}_{\alpha\in A}$ be an
open covering of $M$ such that $\varphi|_{U_\alpha}=\dd\sigma_\alpha$
for smooth functions $\sigma_\alpha:U_\alpha\to\R$, and define
\begin{equation}
    \Omega_\alpha
    :=
    e^{-\sigma_\alpha}\Omega|_{U_\alpha},
    \qquad\text{equivalently}\qquad
    \Omega|_{U_\alpha}
    =
    e^{\sigma_\alpha}\Omega_\alpha.
    \label{eq:local-symplectic-form}
\end{equation}
Then
$\dd\Omega_\alpha=e^{-\sigma_\alpha}(-\dd\sigma_\alpha\wedge\Omega+\dd\Omega)=0$
by \eqref{eq:LCS-definition}, and $\Omega_\alpha$ is nondegenerate
since $\Omega$ is, so that $(U_\alpha,\Omega_\alpha)$ is a symplectic
manifold. If the Lee form is globally exact, $\varphi=\dd\sigma$, then
$e^{-\sigma}\Omega$ is a globally defined symplectic form and the LCS
structure is said to be \emph{globally conformally symplectic}. The
distinction between local and global exactness of the Lee form is one
of the geometric reasons why compatibility between local descriptions
is important.

\paragraph{Overlap relations and the conformal cocycle.}
On each connected component of an overlap $U_\alpha\cap U_\beta$, the
local representatives of the globally defined Lee form agree,
$\dd\sigma_\alpha=\dd\sigma_\beta$, so that
\begin{equation}
    \sigma_\alpha-\sigma_\beta
    =
    C_{\alpha\beta},
    \qquad
    c_{\alpha\beta}
    :=
    e^{C_{\alpha\beta}}>0,
    \qquad
    \Omega_\beta
    =
    c_{\alpha\beta}\Omega_\alpha,
    \label{eq:symplectic-overlap}
\end{equation}
where $C_{\alpha\beta}\in\R$ is constant and the last identity follows
from $e^{\sigma_\alpha}\Omega_\alpha=e^{\sigma_\beta}\Omega_\beta$ on
the overlap. On a triple overlap the constants satisfy
$C_{\alpha\beta}+C_{\beta\gamma}=C_{\alpha\gamma}$, which implies
\begin{equation}
    c_{\alpha\beta}c_{\beta\gamma}
    =
    c_{\alpha\gamma}.
    \label{eq:conformal-cocycle}
\end{equation}
Thus, the constants $c_{\alpha\beta}$ define a multiplicative
conformal cocycle relating the local symplectic representatives.

\paragraph{LCS Hamiltonian dynamics.}
Consider a smooth Hamiltonian function $H:M\to\R$. The LCS Hamiltonian
vector field $\xi_H$ is defined intrinsically by
\begin{equation}
    \iota_{\xi_H}\Omega
    =
    \dd_\varphi H
    =
    \dd H-H\varphi,
    \label{eq:LCS-Hamiltonian-vector-field}
\end{equation}
which determines $\xi_H$ uniquely since $\Omega$ is nondegenerate. On
the local chart $U_\alpha$, define the local Hamiltonian function
\begin{equation}
    H_\alpha
    :=
    e^{-\sigma_\alpha}H|_{U_\alpha}.
    \label{eq:local-Hamiltonian}
\end{equation}
Substituting $\Omega=e^{\sigma_\alpha}\Omega_\alpha$,
$\varphi=\dd\sigma_\alpha$ and $H=e^{\sigma_\alpha}H_\alpha$ into
\eqref{eq:LCS-Hamiltonian-vector-field}, the terms containing
$\dd\sigma_\alpha$ cancel and one obtains
$\iota_{\xi_H}\Omega_\alpha=\dd H_\alpha$. Therefore, the restriction
of $\xi_H$ to $U_\alpha$ is the ordinary Hamiltonian vector field
$X_{H_\alpha}$ of $H_\alpha$ with respect to $\Omega_\alpha$: on each local
symplectic chart, the LCS Hamiltonian dynamics is represented by an
ordinary Hamiltonian system.

\subsection{Local canonical representation}
\label{subsec:local-canonical}

For the Hamiltonian variational construction developed in this paper,
we consider local cotangent-type symplectic representations \cite{Marsden1999}. Throughout
the paper we take $M=T^*Q$ with bundle projection $\pi_Q:T^*Q\to Q$,
denote the LCS variables by $(q^i,p_i)$,
write $\Theta_{\mathrm{can}}=p_i\,\dd q^i$ for the Liouville one-form,
assume that the Lee form is \emph{basic}, i.e.,
\begin{equation}
    \varphi
    =
    \pi_Q^*\varphi_Q
    \label{eq:basic-Lee-form}
\end{equation}
for a closed one-form $\varphi_Q$ on $Q$, so that the local potentials
depend on the base point only, $\sigma_\alpha=\sigma_\alpha(q)$, and
assume that the LCS two-form is of \emph{cotangent type} \cite{HaRy99}, namely
\begin{equation}
    \Omega
    =
    -\dd_\varphi\Theta_{\mathrm{can}}
    =
    \dd q^i\wedge\dd p_i
    +
    \varphi\wedge\Theta_{\mathrm{can}} .
    \label{eq:cotangent-type-hypothesis}
\end{equation}
In other words, $(T^*Q,-\dd_\varphi\Theta_{\mathrm{can}})$ is the
cotangent bundle twisted by the basic closed one-form $\varphi$
\cite{HaRy99,vaisman}.
These are standing hypotheses, not consequences of
\eqref{eq:LCS-definition}, and they are what makes the local
representation below canonical. Note that $\dd_\varphi\Omega
=-\dd_\varphi^2\Theta_{\mathrm{can}}=0$, so
\eqref{eq:cotangent-type-hypothesis} does define an LCS structure for
every closed basic $\varphi$; nondegeneracy holds because, $\varphi$
being basic, the additional term
$\varphi\wedge\Theta_{\mathrm{can}}=\varphi_jp_i\,\dd q^j\wedge\dd q^i$
is a two-form in the base differentials only, so that the matrix of
$\Omega$ in the coordinates $(q^i,p_i)$ is
$\begin{pmatrix}A&I\\-I&0\end{pmatrix}$ with $A$ antisymmetric, which
has determinant one. The restriction to basic Lee forms is essential
here: for a closed one-form with fiber components, such as
$\varphi=\dd p$ on $T^*\R$, the two-form
\eqref{eq:cotangent-type-hypothesis} equals $(1-p)\,\dd q\wedge\dd p$
and degenerates along $\{p=1\}$.

We further assume that the local representations are adapted to the
bundle structure: the chart domains are of the form
$U_\alpha=\pi_Q^{-1}(V_\alpha)$, where $\{V_\alpha\}$ is an open
covering of $Q$ by domains of the local Lee potentials,
$\varphi_Q|_{V_\alpha}=\dd\sigma_\alpha$. For the initial construction, we assume that on overlaps the base
coordinates $q^i$ coincide. Under this assumption, the transition maps
between local canonical representations only rescale the fiber
variables. Under a nonlinear change of base coordinates, the transition
maps would be nonlinear cotangent lifts, to which the Runge--Kutta and
Galerkin results below do not apply as stated; we extend the construction to affine base transitions, which arise for instance from angular coordinates, in Section~\ref{subsec:affine}. 

In such a representation, the local symplectic form on $U_\alpha$ is
written in canonical coordinates $(q^i,r_i)$ as
$\Omega_\alpha=\dd q^i\wedge\dd r_i$, the local Hamilton equations
associated with $H_\alpha$ are $\dot q^i=\partial H_\alpha/\partial r_i$
and $\dot r_i=-\partial H_\alpha/\partial q^i$, and the local canonical
momenta are related to the global momenta by
\begin{equation}
    r_i
    =
    e^{-\sigma_\alpha(q)}p_i,
    \qquad
    \Psi_\alpha(q,p)
    :=
    \left(
        q,
        e^{-\sigma_\alpha(q)}p
    \right),
    \label{eq:local-global-momentum}
\end{equation}
where $\Psi_\alpha$ denotes this local change of representation. Under
the hypothesis \eqref{eq:cotangent-type-hypothesis} the two
descriptions are consistent: with $\varphi=\dd\sigma_\alpha$ on
$U_\alpha$,
\begin{equation}
    \dd q^i\wedge\dd r_i
    =
    e^{-\sigma_\alpha}
    \left(
        \dd q^i\wedge\dd p_i
        -
        p_i\,\dd q^i\wedge\dd\sigma_\alpha
    \right)
    =
    e^{-\sigma_\alpha}
    \left(
        \dd q^i\wedge\dd p_i
        +
        \varphi\wedge\Theta_{\mathrm{can}}
    \right)
    =
    e^{-\sigma_\alpha}\Omega,
    \label{eq:cotangent-consistency}
\end{equation}
which is \eqref{eq:local-symplectic-form} expressed in the canonical
coordinates $(q^i,r_i)$; equivalently,
$\Psi_\alpha^*\Omega_\alpha=e^{-\sigma_\alpha}\Omega$.

Consider now the local coordinates $(q^i,r_i)$ on $U_\alpha$ and
$(q^i,s_i)$ on $U_\beta$. On the overlap, $r=e^{-\sigma_\alpha}p$ and
$s=e^{-\sigma_\beta}p$, so that $s=c_{\alpha\beta}r$ by
\eqref{eq:symplectic-overlap}, and the transition between the two local
canonical representations is
\begin{equation}
    T_{\alpha\beta}
    :=
    \Psi_\beta\circ\Psi_\alpha^{-1},
    \qquad
    T_{\alpha\beta}(q,r)
    =
    (q,c_{\alpha\beta}r),
    \qquad
    T_{\alpha\beta}^*\Omega_\beta
    =
    c_{\alpha\beta}\Omega_\alpha.
    \label{eq:T-alpha-beta}
\end{equation}
Thus $T_{\alpha\beta}$ is conformally symplectic rather than
symplectic unless $c_{\alpha\beta}=1$. The local Hamiltonians satisfy a
corresponding transformation law on the overlap: since $H$ is globally
defined, $e^{\sigma_\alpha}H_\alpha=e^{\sigma_\beta}H_\beta$, and
taking into account the change of momentum variable
$s=c_{\alpha\beta}r$ we obtain
\begin{equation}
    H_\beta(q,s)
    =
    c_{\alpha\beta}
    H_\alpha
    \left(
        q,c_{\alpha\beta}^{-1}s
    \right),
    \label{eq:local-Hamiltonian-overlap}
\end{equation}
and, differentiating,
\begin{equation}
    D_1H_\beta(q,s)
    =
    c_{\alpha\beta}
    D_1H_\alpha
    \left(
        q,c_{\alpha\beta}^{-1}s
    \right),
    \qquad
    D_2H_\beta(q,s)
    =
    D_2H_\alpha
    \left(
        q,c_{\alpha\beta}^{-1}s
    \right).
    \label{eq:H-overlap-D1}
\end{equation}

\paragraph{Equivariance of the continuous Hamiltonian dynamics.}
The local Hamiltonian vector fields $X_{H_\alpha}$ and $X_{H_\beta}$
are related by the transition map. Indeed, by \eqref{eq:T-alpha-beta}
and \eqref{eq:local-Hamiltonian-overlap},
\[
    T_{\alpha\beta}^*
    \big(
        \iota_{T_{\alpha\beta*}X_{H_\alpha}}\Omega_\beta
    \big)
    =
    \iota_{X_{H_\alpha}}T_{\alpha\beta}^*\Omega_\beta
    =
    c_{\alpha\beta}\,\dd H_\alpha
    =
    \dd(H_\beta\circ T_{\alpha\beta})
    =
    T_{\alpha\beta}^*\dd H_\beta,
\]
so that $\iota_{T_{\alpha\beta*}X_{H_\alpha}}\Omega_\beta=\dd H_\beta$,
i.e.,
\begin{equation}
    X_{H_\beta}\circ T_{\alpha\beta}
    =
    \mathrm{T}T_{\alpha\beta}\circ X_{H_\alpha};
    \label{eq:continuous-vector-field-equivariance}
\end{equation}
in coordinates, both sides equal
$\big(D_2H_\alpha(q,r),-c_{\alpha\beta}D_1H_\alpha(q,r)\big)$ by
\eqref{eq:H-overlap-D1}. If $\Phi_{\alpha,t}$ and $\Phi_{\beta,t}$
denote the local Hamiltonian flows of $X_{H_\alpha}$ and
$X_{H_\beta}$, uniqueness of solutions implies
\begin{equation}
    \Phi_{\beta,t}\circ T_{\alpha\beta}
    =
    T_{\alpha\beta}\circ\Phi_{\alpha,t}.
    \label{eq:continuous-flow-equivariance}
\end{equation}
Thus, the two local Hamiltonian descriptions are related by the
transition map and represent the same dynamics on the overlap.

\subsection{Type II Hamiltonian variational principle}

We now recall the Type~II Hamiltonian variational principle in the
local canonical representation, following \cite{LeZh11}. This
formulation is particularly appropriate for the discrete construction
because its discrete counterpart is expressed in terms of right
discrete Hamiltonians. Fix a local symplectic representative
$(U_\alpha,\Omega_\alpha)$ with canonical coordinates $(q,r)$ and local
Hamiltonian $H_\alpha(q,r)$. For curves $(q(t),r(t))$, $t\in[0,T]$,
satisfying the mixed boundary conditions $q(0)=q_0$ and $r(T)=r_T$, we
define the action by
\begin{equation}
    \mathcal S_\alpha^{II}[q,r]
    =
    r(T)\cdot q(T)
    -
    \int_0^T
    \left[
        r(t)\cdot\dot q(t)
        -
        H_\alpha(q(t),r(t))
    \right]
    \dd t.
    \label{eq:continuous-Type-II-action}
\end{equation}

\begin{proposition}
\label{prop:continuous-Type-II}
A curve $(q(t),r(t))$ is stationary for
\eqref{eq:continuous-Type-II-action} with respect to variations
satisfying $\delta q(0)=0$ and $\delta r(T)=0$ if and only if it
satisfies Hamilton's equations
\begin{equation}\label{loc-Ham-eqn}
    \dot{q}(t)
    =
    \frac{\partial H_\alpha}{\partial r}
    (q(t),r(t)),
    \qquad
    \dot{r}(t)
    =
    -\frac{\partial H_\alpha}{\partial q}
    (q(t),r(t)).
\end{equation}
\end{proposition}

\begin{proof}
Integrating the term $r\cdot\delta\dot q$ by parts and using
$\delta q(0)=0$ and $\delta r(T)=0$, the boundary contributions cancel
and
\[
    \delta\mathcal S_\alpha^{II}
    =
    \int_0^T
    \left[
        \left(
            \dot r
            +
            \frac{\partial H_\alpha}{\partial q}
        \right)
        \cdot\delta q
        -
        \left(
            \dot q
            -
            \frac{\partial H_\alpha}{\partial r}
        \right)
        \cdot\delta r
    \right]
    \dd t,
\]
which vanishes for arbitrary interior variations if and only if
\eqref{loc-Ham-eqn} holds.
\end{proof}

\noindent
This is the standard Hamiltonian variational principle with Type~II
boundary conditions. Its importance for the present work is that it
uses the initial position and final momentum as boundary data, which
are precisely the arguments of a right Hamiltonian generating
function.

\paragraph{Recovery of the LCS Hamilton equations.}
The local canonical Hamilton equations can be expressed in the original
LCS variables $(q,p)$ through $H_\alpha(q,r)=e^{-\sigma_\alpha(q)}H(q,p)$
and $r_i=e^{-\sigma_\alpha(q)}p_i$. Writing
$\varphi=\varphi_i(q)\,\dd q^i$ with
$\varphi_i=\partial\sigma_\alpha/\partial q^i$, differentiation of
$r_i=e^{-\sigma_\alpha}p_i$ gives
$\dot r_i=e^{-\sigma_\alpha}(\dot p_i-\varphi_j\dot q^jp_i)$, and
substituting into \eqref{loc-Ham-eqn} we obtain
\begin{equation}
\begin{split}
    \dot q^i
    &=
    \frac{\partial H}{\partial p_i},
    \\
    \dot p_i
    &=
    -
    \frac{\partial H}{\partial q^i}
    -
    A_{ij}
    \frac{\partial H}{\partial p_j}
    +
    H\varphi_i,
    \qquad
    A_{ij}
    :=
    \varphi_i p_j-\varphi_jp_i.
    \label{eq:LCS-q-p-equation}
\end{split}
\end{equation}
These equations give the local coordinate expression of the intrinsic
LCS Hamiltonian equation \eqref{eq:LCS-Hamiltonian-vector-field}.
Hence, the Type~II variational principle does not introduce a different
continuous dynamics; rather, it provides a variational formulation of
the LCS Hamiltonian dynamics, previously described in \cite{EsGeGu}, by
expressing the local dynamics in canonical Hamiltonian form.

\paragraph{Conformal preservation of the continuous flow.}
A characteristic property of LCS Hamiltonian dynamics is that its
Hamiltonian flow preserves the LCS two-form up to a conformal factor.
By Cartan's formula, \eqref{eq:LCS-definition} and
\eqref{eq:LCS-Hamiltonian-vector-field},
\begin{equation}
    \Lie_{\xi_H}\Omega
    =
    \dd(\dd H-H\varphi)
    +
    \iota_{\xi_H}(\varphi\wedge\Omega)
    =
    -\dd H\wedge\varphi
    +
    \varphi(\xi_H)\Omega
    -
    \varphi\wedge(\dd H-H\varphi)
    =
    \varphi(\xi_H)\Omega,
    \label{eq:LCS-Lie-preservation}
\end{equation}
where we have used $\dd\varphi=0$ and $\varphi\wedge\varphi=0$. Let
$\psi_t:M\to M$ denote the integral flow of $\xi_H$. Since
$\frac{\dd}{\dd t}\psi_t^*\Omega=\psi_t^*(\Lie_{\xi_H}\Omega)
=(\varphi(\xi_H)\circ\psi_t)\,\psi_t^*\Omega$ and $\psi_0=\id$,
\begin{equation}
    \psi_t^*\Omega
    =
    \exp\left(
        \int_0^t
        \varphi(\xi_H)(\psi_s(x))\,\dd s
    \right)
    \Omega.
    \label{eq:continuous-conformal-factor-general}
\end{equation}
If the flow segment under consideration remains in a local chart
$U_\alpha$, then $\varphi(\xi_H)=\dd\sigma_\alpha(\xi_H)$, so that
$\varphi(\xi_H)(\psi_t(x))$ is the time derivative of
$\sigma_\alpha(\psi_t(x))$, and therefore
\begin{equation}
    \psi_t^*\Omega
    =
    e^{
        \sigma_\alpha\circ\psi_t
        -
        \sigma_\alpha
    }
    \Omega.
    \label{eq:continuous-LCS-preservation}
\end{equation}
Although the right-hand side is expressed in terms of the local
potential $\sigma_\alpha$, the conformal factor is independent of the
choice of local potential: whenever $x$ and $\psi_t(x)$ belong to the
same connected component of $U_\alpha\cap U_\beta$, the constant
$C_{\alpha\beta}$ cancels and
$\sigma_\alpha(\psi_t(x))-\sigma_\alpha(x)
=\sigma_\beta(\psi_t(x))-\sigma_\beta(x)$. Equation
\eqref{eq:continuous-LCS-preservation} is the continuous
structure-preservation law that the discrete theory will ultimately
seek to reproduce.

\subsection{Continuous example}

We conclude the continuous discussion with the example that will later
be used to illustrate the numerical constructions.

\begin{example}
\label{ex:continuous-LCS-oscillator}
Let $Q=\R$ and take $Q\times\R$ as the base manifold, that is, the
setting of Section~\ref{subsec:local-canonical} with $Q$ replaced by
$Q\times\R$, with coordinates $(q,x)$ and induced coordinates
$(q,p,x,y)$ on the cotangent bundle $T^*(Q\times\R)$. Consider the Lee
form $\varphi=-\dd x$, with local conformal potential $\sigma=-x$.
Since $\varphi$ is the pullback of a closed one-form on $Q\times\R$
under the cotangent-bundle projection, the cotangent-type hypothesis
\eqref{eq:cotangent-type-hypothesis} applies, and in the LCS variables
\begin{equation}
    \Omega
    =
    \dd q\wedge\dd p
    +
    \dd x\wedge\dd y
    +
    p\,\dd q\wedge\dd x .
    \label{eq:example-Omega-LCS-variables}
\end{equation}
Consider the Hamiltonian
$H(q,p,x,y)=\tfrac12(p^2+y^2)+\tfrac{k}{2}q^2+\mu x$, where $k$ and
$\mu$ are real parameters. Introducing the
local canonical momenta $r=e^xp$ and $s=e^xy$, the local symplectic
form is $\Omega_{\mathrm{loc}}=\dd q\wedge\dd r+\dd x\wedge\dd s$, with
$\Omega=e^{-x}\Omega_{\mathrm{loc}}$ consistently with
\eqref{eq:example-Omega-LCS-variables}, and the local Hamiltonian is
\begin{equation}
    H_{\mathrm{loc}}(q,r,x,s)
    =
    \frac12e^{-x}(r^2+s^2)
    +
    e^x
    \left(
        \frac{k}{2}q^2+\mu x
    \right),
    \label{eq:example-local-H}
\end{equation}
whose canonical Hamilton equations are
\begin{equation}
    \dot q
    =
    e^{-x}r,
    \qquad
    \dot r
    =
    -ke^x q,
    \qquad
    \dot x
    =
    e^{-x}s,
    \qquad
    \dot s
    =
    \frac12e^{-x}(r^2+s^2)
    -
    e^x
    \left[
        \frac{k}{2}q^2+\mu(x+1)
    \right].
    \label{eq:example-continuous-Hamilton}
\end{equation}
The continuous LCS flow therefore satisfies
$\psi_t^*\Omega=e^{-(x(t)-x(0))}\Omega$. Its numerical discretization
will be considered only after the general discrete theory has been
established.
\end{example}

\noindent
The continuous framework therefore supplies two pieces of information
that must survive discretization: each local representative is an
ordinary symplectic Hamiltonian system, and different representatives
are related on overlaps by the constant fiber scalings
$T_{\alpha\beta}$. The exact flows respect both structures
automatically. We now construct Type~II discrete Hamiltonian maps
locally and then determine the additional compatibility condition
required for those maps to globalize.


\section{Discrete Type II Hamiltonian mechanics}
\label{sec:discrete-ham}

We now turn to the discrete problem. The construction is first carried
out on a single local symplectic representative
$(U_\alpha,\Omega_\alpha)$, with canonical coordinates $(q,r)$ and
local Hamiltonian $H_\alpha$; no compatibility between different
charts is imposed at this stage. The basic object is a right discrete
Hamiltonian, that is, a Type~II generating function depending on an
initial position and a final momentum. Such generating-function
techniques have also been used to construct symplectic integrators
for optimal control problems \cite{LeDiSa-07}. We first recall the exact
discrete right Hamiltonian and then replace it by a computable
approximation, which produces a local symplectic one-step map. The
behavior of these local constructions on overlaps will be studied in
Section~\ref{sec:globalization}.

\subsection{The exact discrete right Hamiltonian}

Fix a time step $h>0$ and consider the exact solution $(q(t),r(t))$,
$t\in[0,h]$, of the local Hamiltonian system $\dot q=D_2H_\alpha(q,r)$,
$\dot r=-D_1H_\alpha(q,r)$, subject to the Type~II boundary conditions
\begin{equation}
    q(0)=q_0,
    \qquad
    r(h)=r_1.
    \label{eq:exact-Type-II-boundary}
\end{equation}
We assume throughout this discussion that the corresponding local
two-point boundary-value problem has a unique solution.

\begin{definition}
The \emph{exact discrete right Hamiltonian} associated with
$H_\alpha$ is
\begin{equation}
    H_{\alpha,h}^{+,E}(q_0,r_1)
    =
    r_1\cdot q_1
    -
    \int_0^h
    \left[
        r(t)\cdot\dot q(t)
        -
        H_\alpha(q(t),r(t))
    \right]
    \dd t,
    \label{eq:exact-discrete-right-Hamiltonian}
\end{equation}
where $q_1=q(h)$ and $(q(t),r(t))$ is the exact solution satisfying
\eqref{eq:exact-Type-II-boundary}.
\end{definition}
\noindent
The exact discrete right Hamiltonian is the Type~II generating function
of the exact time-$h$ Hamiltonian flow.

\begin{proposition}
\label{prop:exact-discrete-right-H}
The exact discrete right Hamiltonian satisfies
\begin{equation}
    r_0
    =
    D_1H_{\alpha,h}^{+,E}(q_0,r_1),
    \qquad
    q_1
    =
    D_2H_{\alpha,h}^{+,E}(q_0,r_1),
    \label{eq:exact-right-discrete-Hamilton-equations}
\end{equation}
where $r_0=r(0)$ and $q_1=q(h)$. Consequently, $H_{\alpha,h}^{+,E}$
generates the exact time-$h$ Hamiltonian map.
\end{proposition}

\begin{proof}
Evaluating the variation of the Type~II action
\eqref{eq:continuous-Type-II-action} on a solution of Hamilton's
equations, the interior terms vanish and only the boundary
contribution
$\delta H_{\alpha,h}^{+,E}=r_0\cdot\delta q_0+q_1\cdot\delta r_1$
remains, which is \eqref{eq:exact-right-discrete-Hamilton-equations}.
\end{proof}
\noindent
Thus the exact time-$h$ Hamiltonian map may be encoded entirely by a
function of the mixed endpoint data $(q_0,r_1)$.

\paragraph{Approximate discrete right Hamiltonians.}
The exact discrete Hamiltonian
\eqref{eq:exact-discrete-right-Hamiltonian} is generally not available
in closed form. A numerical method is obtained by replacing it with an
approximation.

\begin{definition}
A \emph{discrete right Hamiltonian} is a smooth function
$H_{\alpha,h}^{d+}(q_0,r_1)$ approximating the exact discrete right
Hamiltonian $H_{\alpha,h}^{+,E}(q_0,r_1)$.
\end{definition}
\noindent
The corresponding right discrete Hamilton equations are
\begin{equation}
    r_0
    =
    D_1H_{\alpha,h}^{d+}(q_0,r_1),
    \qquad
    q_1
    =
    D_2H_{\alpha,h}^{d+}(q_0,r_1).
    \label{eq:right-discrete-Hamilton-equations}
\end{equation}
Given $(q_0,r_0)$, the first equation is solved implicitly for $r_1$
and the second equation then determines $q_1$. Whenever this procedure
is locally well posed, it defines a one-step map
$F_{\alpha,h}:(q_0,r_0)\mapsto(q_1,r_1)$. A convenient sufficient
regularity condition is the nonsingularity of the mixed derivative
$D_{12}H_{\alpha,h}^{d+}(q_0,r_1)$, under which the first equation can
be solved locally for $r_1$ as a function of $(q_0,r_0)$ by the
implicit function theorem.

\subsection{Discrete Type II variational formulation}

The one-step relations \eqref{eq:right-discrete-Hamilton-equations}
admit a discrete phase-space variational interpretation. Consider a
sequence $\{(q_k,r_k)\}_{k=0}^N$ with the Type~II endpoint conditions
$q_0$ and $r_N$ fixed, and define the discrete Type~II action by
\begin{equation}
    \mathcal S_{\alpha,d}^{II}
    =
    r_N\cdot q_N
    -
    \sum_{k=0}^{N-1}
    \left[
        r_{k+1}\cdot q_{k+1}
        -
        H_{\alpha,h}^{d+}(q_k,r_{k+1})
    \right].
    \label{eq:discrete-Type-II-action}
\end{equation}
The endpoint term cancels the final contribution $r_N\cdot q_N$
contained in the sum. Taking the variation and using $\delta q_0=0$,
$\delta r_N=0$ gives
\begin{equation}
    \delta\mathcal S_{\alpha,d}^{II}
    =
    \sum_{k=1}^{N-1}
    \left[
        -r_k
        +
        D_1H_{\alpha,h}^{d+}(q_k,r_{k+1})
    \right]\cdot\delta q_k
    +
    \sum_{k=1}^{N-1}
    \left[
        -q_k
        +
        D_2H_{\alpha,h}^{d+}(q_{k-1},r_k)
    \right]\cdot\delta r_k,
    \label{eq:variation-discrete-Type-II}
\end{equation}
and we obtain the following interior stationarity equations.

\begin{proposition}
\label{prop:discrete-Type-II-stationarity}
A sequence is stationary for \eqref{eq:discrete-Type-II-action}, with
$q_0$ and $r_N$ fixed, if and only if its interior variables satisfy
\begin{equation}
    r_k
    =
    D_1H_{\alpha,h}^{d+}(q_k,r_{k+1}),
    \qquad
    q_k
    =
    D_2H_{\alpha,h}^{d+}(q_{k-1},r_k),
    \qquad
    k=1,\ldots,N-1.
    \label{eq:discrete-Type-II-interior-r}
\end{equation}
\end{proposition}

\begin{remark}
\label{rem:boundary-momentum}
The stationarity equations are interior equations. Because $q_0$ is
fixed, variation of the discrete action does not yield an equation for
the initial momentum $r_0$, and because $q_N$ cancels from
\eqref{eq:discrete-Type-II-action}, it does not yield an equation for
the final position either. Both are recovered from the Type~II
generating-function relations $r_0=D_1H_{\alpha,h}^{d+}(q_0,r_1)$ and
$q_N=D_2H_{\alpha,h}^{d+}(q_{N-1},r_N)$. Likewise, the one-step update
is described by the pair of generating relations
\begin{equation}\label{eq:loc-dHEq}
     r_k
    =
    D_1H_{\alpha,h}^{d+}(q_k,r_{k+1}),
    \qquad
    q_{k+1}
    =
    D_2H_{\alpha,h}^{d+}(q_k,r_{k+1}).
\end{equation}
This distinction between the interior stationarity conditions and the
endpoint generating relations is useful when formulating the discrete
Type~II principle with mixed boundary data.
\end{remark}

\paragraph{Local symplecticity.}

\begin{proposition}
\label{prop:local-discrete-symplecticity}
Suppose that $H_{\alpha,h}^{d+}$ is regular enough to define the local
one-step map $F_{\alpha,h}$. Then $F_{\alpha,h}$ is symplectic with
respect to the local symplectic form $\Omega_\alpha$.
\end{proposition}

\begin{proof}
By \eqref{eq:right-discrete-Hamilton-equations},
$\dd H_{\alpha,h}^{d+}=r_0\,\dd q_0+q_1\,\dd r_1$, and taking one more
exterior derivative gives
$0=\dd r_0\wedge\dd q_0+\dd q_1\wedge\dd r_1$, i.e.,
\begin{equation}
    \dd q_0\wedge\dd r_0
    =
    \dd q_1\wedge\dd r_1,
    \label{eq:local-symplecticity-coordinate}
\end{equation}
which is $F_{\alpha,h}^*\Omega_\alpha=\Omega_\alpha$ for
$\Omega_\alpha=\dd q^i\wedge\dd r_i$.
\end{proof}
\noindent
This is a purely local statement. It says that the discrete dynamics is
symplectic when expressed in one canonical representative. It does not
yet imply that independently constructed maps on different local
representations determine a single global numerical map.

\paragraph{Order of a discrete Hamiltonian method.}

The approximation order of the discrete Hamiltonian controls the order
of the induced one-step map.

\begin{theorem}[\cite{ScLe18}]
\label{thm:discrete-Hamiltonian-order}
Let $H_{\alpha,h}^{+,E}$ denote the exact discrete right Hamiltonian and
let $H_{\alpha,h}^{d+}$ be an approximate discrete right Hamiltonian.
Assume, on a compact set on which the exact and discrete maps are
well-defined, that
\begin{equation}
    H_{\alpha,h}^{d+}
    -
    H_{\alpha,h}^{+,E}
    =
    \cO(h^{p+1})
    \label{eq:Hd-order}
\end{equation}
in the $C^1$ sense as $h\to0$. Then, under the usual local regularity
assumptions for the corresponding generating relations, the one-step
map generated by $H_{\alpha,h}^{d+}$ approximates the exact
time-$h$ Hamiltonian flow with order $p$.
\end{theorem}

\noindent
The $C^1$ requirement is relevant because the numerical map is
generated by the derivatives
\[
    D_1H_{\alpha,h}^{d+},
    \qquad
    D_2H_{\alpha,h}^{d+}
\]
rather than by the value of the discrete Hamiltonian alone.

\begin{remark}[Converse of Theorem~\ref{thm:discrete-Hamiltonian-order}]
\label{rem:order-converse}
The converse also holds, up to an additive constant. Suppose that the
one-step map $F_{\alpha,h}$ generated by a regular discrete right
Hamiltonian $H^{d+}_{\alpha,h}$ approximates the exact time-$h$ flow
$\Phi_{\alpha,h}$ with order $p$ on a compact set, i.e.,
$F_{\alpha,h}-\Phi_{\alpha,h}=\cO(h^{p+1})$ there.
Let $\pi_q$ and $\pi_r$ denote the projections onto the position
and momentum, respectively, and let
\[
    G_{\alpha,h}(q_0,r_0)
    =
    \bigl(q_0,\pi_rF_{\alpha,h}(q_0,r_0)\bigr)
    =
    (q_0,r_1)
\]
denote the map from the initial data to the mixed data.
Let $G^E_{\alpha,h}$ denote the corresponding mixed-data map
for $\Phi_{\alpha,h}$.
Assume that these maps admit inverse branches on a common
mixed-data neighborhood, with uniformly bounded Lipschitz
constants for the inverses as $h\to0$.
Assume also that both inverse branches take values in the
compact set on which
$F_{\alpha,h}-\Phi_{\alpha,h}=\cO(h^{p+1})$
holds uniformly. The generating relations
\eqref{eq:exact-right-discrete-Hamilton-equations} and
\eqref{eq:right-discrete-Hamilton-equations} can then be written as
\[
    D_1H^{d+}_{\alpha,h}
    =
    \pi_r\circ G_{\alpha,h}^{-1},
    \qquad
    D_2H^{d+}_{\alpha,h}
    =
    \pi_q\circ F_{\alpha,h}\circ G_{\alpha,h}^{-1},
\]
and likewise for $H^{+,E}_{\alpha,h}$ with $\Phi_{\alpha,h}$ and
$G^E_{\alpha,h}$ in place of $F_{\alpha,h}$ and $G_{\alpha,h}$.
Since $G_{\alpha,h}-G^E_{\alpha,h}=\cO(h^{p+1})$,
the inverses and the compositions differ by $\cO(h^{p+1})$
as well \cite[Lemma~2.3]{ScLe18}, so that
\[
    D_iH^{d+}_{\alpha,h}
    -
    D_iH^{+,E}_{\alpha,h}
    =
    \cO(h^{p+1}),
    \qquad i=1,2,
\]
uniformly on compact subsets of the common mixed-data
neighborhood. On a fixed convex compact subdomain of this
neighborhood, integrating the derivative difference along
line segments shows that, up to an additive constant
depending only on $h$,
$H^{d+}_{\alpha,h}-H^{+,E}_{\alpha,h}
=\cO(h^{p+1})$ in the $C^1$ sense.
Since additive constants do not affect the generated map,
the order of a Hamiltonian variational integrator and the
order to which its discrete right Hamiltonian approximates
the exact discrete right Hamiltonian coincide.
\end{remark}

\paragraph{The local construction and the compatibility problem.}

Thus, on every local symplectic representative
$(U_\alpha,\Omega_\alpha,H_\alpha)$, an approximate discrete right
Hamiltonian $H_{\alpha,h}^{d+}$ determines, under the regularity
assumptions above, a symplectic one-step map $F_{\alpha,h}$. On an
overlap, the continuous Hamiltonians are related by
\eqref{eq:local-Hamiltonian-overlap} and their exact flows are
equivariant with respect to $T_{\alpha\beta}$. There is no analogous
automatic conclusion for two independently chosen approximations
$H_{\alpha,h}^{d+}$ and $H_{\beta,h}^{d+}$; in particular, local
symplecticity does not imply
$F_{\beta,h}\circ T_{\alpha\beta}=T_{\alpha\beta}\circ F_{\alpha,h}$.
The missing ingredient is a compatibility condition on the discrete
generating functions. We derive it next from the transformation law of
the exact discrete right Hamiltonian; this keeps the discrete
condition tied directly to the underlying continuous LCS geometry.


\section{Globalization of discrete Hamiltonian dynamics}
\label{sec:globalization}

The construction of Section~\ref{sec:discrete-ham} associates a
symplectic one-step map with each local canonical representation of the
LCS Hamiltonian system. The remaining question is whether these local
maps represent a single discrete dynamics on the underlying LCS phase
space. The continuous theory provides the model for the required
compatibility: on a connected component of an overlap
$U_\alpha\cap U_\beta$, the local canonical variables are related by
the constant fiber scaling $T_{\alpha\beta}$ of \eqref{eq:T-alpha-beta}
and the local Hamiltonians by \eqref{eq:local-Hamiltonian-overlap}. We
first show that the exact discrete right Hamiltonians inherit precisely
the same transformation law. This observation provides the natural
compatibility condition for approximate discrete Hamiltonians.

\subsection{Transformation of the exact discrete right Hamiltonian}

Let $H_{\alpha,h}^{+,E}(q_0,r_1)$ and $H_{\beta,h}^{+,E}(q_0,s_1)$
denote the exact discrete right Hamiltonians
\eqref{eq:exact-discrete-right-Hamiltonian} associated with $H_\alpha$
and $H_\beta$, respectively.

\begin{theorem}
\label{thm:exact-discrete-H-transformation}
For mixed data $(q_0,r_1)$ such that the solution $(q(t),r(t))$
of the $\alpha$-Hamiltonian system with $q(0)=q_0$ and
$r(h)=r_1$ remains in a single connected component of
$U_\alpha\cap U_\beta$ for all $t\in[0,h]$, set
$s_1=c_{\alpha\beta}r_1$. Then the exact discrete
right Hamiltonians satisfy
\begin{equation}
    H_{\beta,h}^{+,E}(q_0,s_1)
    =
    c_{\alpha\beta}
    H_{\alpha,h}^{+,E}
    \left(
        q_0,c_{\alpha\beta}^{-1}s_1
    \right).
    \label{eq:exact-discrete-H-transformation}
\end{equation}
\end{theorem}

\begin{proof}
Let $(q(t),r(t))$ be the solution of the $\alpha$-Hamiltonian system
satisfying $q(0)=q_0$, $r(h)=r_1$, and set $s(t)=c_{\alpha\beta}r(t)$.
By the continuous equivariance relation
\eqref{eq:continuous-flow-equivariance}, $(q(t),s(t))$ is the solution
of the $\beta$-Hamiltonian system with $q(0)=q_0$ and
$s(h)=s_1=c_{\alpha\beta}r_1$. Along these trajectories
$H_\beta(q,s)=c_{\alpha\beta}H_\alpha(q,r)$ by
\eqref{eq:local-Hamiltonian-overlap} and
$s\cdot\dot q=c_{\alpha\beta}\,r\cdot\dot q$, so every term of
\eqref{eq:exact-discrete-right-Hamiltonian} in the
$\beta$-representation is $c_{\alpha\beta}$ times the corresponding
term in the $\alpha$-representation, i.e.,
$H_{\beta,h}^{+,E}(q_0,s_1)=c_{\alpha\beta}H_{\alpha,h}^{+,E}(q_0,r_1)$
with $r_1=c_{\alpha\beta}^{-1}s_1$.
\end{proof}

\noindent
Equation~\eqref{eq:exact-discrete-H-transformation} is not an
additional compatibility assumption. It is inherited by the exact
Type~II generating functions from the overlap relation of the
continuous Hamiltonian systems, and it is the basis for the
compatibility condition imposed below on numerical approximations.

\paragraph{Compatible discrete right Hamiltonians.}

\begin{definition}
\label{def:LCS-compatible-discrete-H}
A collection of local discrete right Hamiltonians
$\{H_{\alpha,h}^{d+}\}_{\alpha\in A}$ is called \emph{LCS-compatible}
if, on every connected component of an overlap $U_\alpha\cap U_\beta$,
one has
\begin{equation}
    H_{\beta,h}^{d+}(q,s)
    =
    c_{\alpha\beta}
    H_{\alpha,h}^{d+}
    \left(
        q,c_{\alpha\beta}^{-1}s
    \right).
    \label{eq:compatible-discrete-H}
\end{equation}
\end{definition}
\noindent
This is the discrete analogue of the continuous transformation law
\eqref{eq:local-Hamiltonian-overlap}. More importantly, it is a
condition imposed directly on the discrete generating functions and can
therefore be checked at the level at which the numerical method is
constructed. Differentiating \eqref{eq:compatible-discrete-H} gives
\begin{equation}
    D_1H_{\beta,h}^{d+}(q,s)
    =
    c_{\alpha\beta}
    D_1H_{\alpha,h}^{d+}
    \left(
        q,c_{\alpha\beta}^{-1}s
    \right),
    \qquad
    D_2H_{\beta,h}^{d+}(q,s)
    =
    D_2H_{\alpha,h}^{d+}
    \left(
        q,c_{\alpha\beta}^{-1}s
    \right);
    \label{eq:compatible-H-D1}
\end{equation}
the cancellation of the factor $c_{\alpha\beta}$ in the second identity
is precisely what is needed for the position update to be independent
of the local representation.

\subsection{Compatibility implies discrete equivariance}

Let $F_{\alpha,h}(q_k,r_k)=(q_{k+1},r_{k+1})$ be the local discrete map
generated by $H_{\alpha,h}^{d+}$ through \eqref{eq:loc-dHEq}, and let
$F_{\beta,h}$ be the map generated by $H_{\beta,h}^{d+}$.

\begin{theorem}[Discrete equivariance]
\label{thm:discrete-equivariance}
Suppose that the local discrete right Hamiltonians are LCS-compatible
in the sense of Definition~\ref{def:LCS-compatible-discrete-H}, and
assume that the corresponding generating relations define unique local
one-step maps. Then
\begin{equation}
    F_{\beta,h}\circ T_{\alpha\beta}
    =
    T_{\alpha\beta}\circ F_{\alpha,h}
    \label{eq:discrete-equivariance}
\end{equation}
on their common domain, i.e., at every point $(q_k,r_k)$ such that
$(q_k,r_k)$ and $F_{\alpha,h}(q_k,r_k)$ lie in the same connected
component of the overlap, so that a single constant $c_{\alpha\beta}$
applies to both endpoints of the step.
\end{theorem}

\begin{proof}
Let $F_{\alpha,h}(q_k,r_k)=(q_{k+1},r_{k+1})$, so that
$r_k=D_1H_{\alpha,h}^{d+}(q_k,r_{k+1})$ and
$q_{k+1}=D_2H_{\alpha,h}^{d+}(q_k,r_{k+1})$, and set
$s_k=c_{\alpha\beta}r_k$, $s_{k+1}=c_{\alpha\beta}r_{k+1}$. By
\eqref{eq:compatible-H-D1},
\[
    D_1H_{\beta,h}^{d+}(q_k,s_{k+1})
    =
    c_{\alpha\beta}D_1H_{\alpha,h}^{d+}(q_k,r_{k+1})
    =
    s_k,
    \qquad
    D_2H_{\beta,h}^{d+}(q_k,s_{k+1})
    =
    D_2H_{\alpha,h}^{d+}(q_k,r_{k+1})
    =
    q_{k+1},
\]
so that $(q_k,s_k)$ and $(q_{k+1},s_{k+1})$ satisfy the generating
relations of $H_{\beta,h}^{d+}$. By uniqueness,
\[
    F_{\beta,h}\big(T_{\alpha\beta}(q_k,r_k)\big)
    =
    (q_{k+1},c_{\alpha\beta}r_{k+1})
    =
    T_{\alpha\beta}\big(F_{\alpha,h}(q_k,r_k)\big),
\]
which proves \eqref{eq:discrete-equivariance}.
\end{proof}
\noindent
Thus the transformation law for the discrete generating functions is a
sufficient condition for compatibility of the resulting numerical
maps.

\paragraph{Gluing of the local discrete maps.}

\begin{lemma}[Gluing lemma]
\label{lem:discrete-gluing}
Let $\{F_{\alpha,h}\}$ be a collection of local maps satisfying
$F_{\beta,h}\circ T_{\alpha\beta}=T_{\alpha\beta}\circ F_{\alpha,h}$
on every overlap. Then the local maps define a unique discrete map on
the phase space covered by the local charts, wherever the local
one-step maps are defined.
\end{lemma}

\begin{proof}
If a point is represented in two overlapping local descriptions by
$z_\beta=T_{\alpha\beta}(z_\alpha)$, then
$F_{\beta,h}(z_\beta)=T_{\alpha\beta}(F_{\alpha,h}(z_\alpha))$, so the
two local updates represent the same point after one time step; the
cocycle relation \eqref{eq:conformal-cocycle} guarantees consistency on
multiple overlaps.
\end{proof}
\noindent
The gluing statement itself is elementary. The relevant point for the
present construction is that Theorem~\ref{thm:discrete-equivariance}
provides a concrete sufficient condition for its hypothesis directly
in terms of the discrete Hamiltonians.

The requirement that each step remain within a single chart
is ensured on compact sets for sufficiently small step sizes;
see Remark~\ref{rem:chart-boundaries}.

\paragraph{Reconstruction in the LCS variables.}

We now return from the local canonical variables to the original LCS
variables through the representation map $\Psi_\alpha$ of
\eqref{eq:local-global-momentum}, whose inverse is
$\Psi_\alpha^{-1}(q,r)=(q,e^{\sigma_\alpha(q)}r)$. Given a local
discrete map $F_{\alpha,h}$ in canonical variables, define the
reconstructed map in the LCS variables by
\begin{equation}
    F_h
    =
    \Psi_\alpha^{-1}
    \circ
    F_{\alpha,h}
    \circ
    \Psi_\alpha.
    \label{eq:reconstructed-global-map}
\end{equation}
Theorem~\ref{thm:discrete-equivariance} guarantees that this definition
is independent of the chosen local representation: on an overlap,
$T_{\alpha\beta}=\Psi_\beta\circ\Psi_\alpha^{-1}$, and therefore
\begin{equation}
    \Psi_\beta^{-1}
    \circ F_{\beta,h}
    \circ\Psi_\beta
    =
    \Psi_\beta^{-1}
    \circ
    T_{\alpha\beta}
    \circ
    F_{\alpha,h}
    \circ
    T_{\alpha\beta}^{-1}
    \circ
    \Psi_\beta
    =
    \Psi_\alpha^{-1}
    \circ
    F_{\alpha,h}
    \circ
    \Psi_\alpha.
\end{equation}
Hence the reconstructed local expressions glue to a well-defined map
$F_h:M\to M$ on the domain reached by the numerical evolution,
characterized on each chart by
\begin{equation}
    \Psi_\alpha\circ F_h
    =
    F_{\alpha,h}\circ\Psi_\alpha.
    \label{eq:commuting-local}
\end{equation}
This identity will be used immediately to transfer local symplecticity
to the conformal preservation law for the global LCS form.

The construction of $F_h$ also admits a symplectic covering
interpretation. The following remark explains this connection
and the role of our chart-wise Type~II formulation.

\begin{remark}[Symplectic covering interpretation]
Consider the cotangent bundle $M=T^*Q$ equipped with the LCS
form $\Omega$ and the basic Lee form
$\varphi=\pi_Q^*\varphi_Q$.
Let $\rho:\widetilde M\to M$ be the minimal covering on which
$\rho^*\varphi=\dd\widetilde\sigma$. Then
\[
    \widetilde\omega=e^{-\widetilde\sigma}\rho^*\Omega,
    \qquad
    \widetilde H=e^{-\widetilde\sigma}\rho^*H
\]
define a symplectic system; see \cite{Otiman-18}.
Since $\widetilde\sigma\circ\gamma-\widetilde\sigma=\lambda_\gamma$ is
constant, equal to a period of $\varphi$, for every deck transformation
$\gamma$, one has $\gamma^*\widetilde\omega=e^{-\lambda_\gamma}
\widetilde\omega$ and $\widetilde H\circ\gamma=e^{-\lambda_\gamma}
\widetilde H$: deck transformations are homotheties rather than
symplectomorphisms, so deck equivariance is not automatic for an
arbitrary symplectic integrator on $\widetilde M$.
A deck transformation is a diffeomorphism
$\gamma:\widetilde M\to\widetilde M$ satisfying
$\rho\circ\gamma=\rho$.
If a symplectic numerical map $\widetilde F_h$ commutes with
every deck transformation, then it defines a map $F_h$ on $M$
satisfying
\[
    \rho\circ\widetilde F_h=F_h\circ\rho.
\]
In Example~\ref{ex:S1-example}, the angular chart transition corresponds
to the deck translation $\theta\mapsto\theta+2\pi$.
Our chart-wise construction is therefore equivalent to a deck-equivariant
symplectic integration on $\widetilde M$. What the Type~II formulation
adds is a verifiable condition for this equivariance at the level of the
discrete generating functions, including the endpoint term under affine
transitions, and a computation that avoids the covering, where the
canonical variables involve factors $e^{\mp\widetilde\sigma}$ that are
unbounded in the number of windings.
\end{remark}

\begin{remark}[Discrete locally conformal Hamiltonian equations]
\label{rem:discrete-LCS-equations}
Following the continuous formulation, we obtain the discrete locally
conformally Hamiltonian dynamics by applying the same local-to-global
construction at the discrete level. Since a discrete right Hamiltonian
depends on the initial position and on the final momentum, the momentum
rescaling \eqref{eq:local-global-momentum} must be evaluated at the
endpoint carrying the momentum. We therefore consider the function
\begin{equation}
    H^{d+}\big|_{U_\alpha}(q_k,p_{k+1};q_{k+1})
    =
    e^{\sigma_\alpha(q_k)}
    H_{\alpha,h}^{d+}
    \big(
        q_k,e^{-\sigma_\alpha(q_{k+1})}p_{k+1}
    \big),
    \label{eq:global-discrete-H}
\end{equation}
in which the final position enters only through the conformal factor
$e^{-\sigma_\alpha(q_{k+1})}$ and is treated as a parameter; the
derivatives $D_1$ and $D_2$ below are taken with respect to the first
two arguments, with $q_{k+1}$ held fixed. Using
\eqref{eq:global-discrete-H} together with
$r_k=e^{-\sigma_\alpha(q_k)} p_k$ in the local discrete Hamiltonian
equations \eqref{eq:loc-dHEq}, we obtain the right discrete locally
conformally Hamiltonian equations
\begin{equation}
\label{r-dLCH}
\begin{cases}
{q}_{k+1}
=
e^{\sigma_{\alpha}({q}_{k+1})
-\sigma_{\alpha}({q}_{k})}
D_2 H^{d+}\big|_{U_\alpha}
({q}_{k},{p}_{k+1};q_{k+1}),
\\[1mm]
{p}_{k}
=
D_1 H^{d+}\big|_{U_\alpha}
({q}_{k},{p}_{k+1};q_{k+1})
-
D\sigma_{\alpha}({q}_{k})
\cdot
H^{d+}\big|_{U_\alpha}
({q}_{k},{p}_{k+1};q_{k+1}).
\end{cases}
\end{equation}
The exponential prefactor in the first equation accounts for the
difference between the momentum rescalings at $q_k$ and at $q_{k+1}$.
The second equation is the discrete counterpart of
$\iota_{\xi_H}\Omega=\dd_\varphi H$, the operator
$D_1-D\sigma_\alpha\,$ acting on the first slot being the discrete
analogue of the Lichnerowicz--de~Rham differential
\eqref{eq:Lichnerowicz-differential}. When $\sigma_\alpha\equiv0$ the
system \eqref{r-dLCH} reduces to \eqref{eq:loc-dHEq}, as it must. Note
that, unlike $H^{d+}_{\alpha,h}$, the function
\eqref{eq:global-discrete-H} is not a Type~II generating function of
the mixed data $(q_k,p_{k+1})$ alone, since the conformal factor at the
final position cannot be expressed in terms of these data; the
variational structure is carried by the local generating function
$H^{d+}_{\alpha,h}$, and \eqref{r-dLCH} is its expression in the LCS
variables.

Thus, the discrete Type~II variational principle yields the
right discrete locally conformally Hamiltonian equations in the
original LCS variables, see \cite{EsGeGu24}.
\end{remark}

\subsection{Global discrete LCS preservation}

We now establish the structure-preservation property of the global
discrete map. Under the cotangent-type hypothesis
\eqref{eq:cotangent-type-hypothesis}, the computation
\eqref{eq:cotangent-consistency} shows that the representation map
$\Psi_\alpha$ satisfies
\begin{equation}
    \Psi_\alpha^*\Omega_\alpha
    =
    e^{-\sigma_\alpha}\Omega,
    \label{eq:Psi-pullback-Omega}
\end{equation}
where the left-hand side is expressed in the LCS variables and
$\Omega_\alpha$ in the local canonical variables.

\begin{theorem}[Global LCS preservation]
\label{thm:global-LCS-preservation}
Let $\{F_{\alpha,h}\}$ be a family of local one-step maps, each
symplectic with respect to $\Omega_\alpha$ and equivariant on overlaps
in the sense of \eqref{eq:discrete-equivariance}, for instance the maps
generated by an LCS-compatible family of regular discrete right
Hamiltonians, and let $F_h$ denote their global reconstruction. Then
\begin{equation}
    F_h^*\Omega
    =
    e^{\sigma_\alpha\circ F_h-\sigma_\alpha}
    \Omega
    \label{eq:global-discrete-LCS-preservation}
\end{equation}
on every chart on which the one-step map is represented.
\end{theorem}

\begin{proof}
By \eqref{eq:Psi-pullback-Omega} and \eqref{eq:commuting-local},
\[
    e^{-\sigma_\alpha\circ F_h}F_h^*\Omega
    =
    F_h^*\big(e^{-\sigma_\alpha}\Omega\big)
    =
    (\Psi_\alpha\circ F_h)^*\Omega_\alpha
    =
    \Psi_\alpha^*F_{\alpha,h}^*\Omega_\alpha
    =
    \Psi_\alpha^*\Omega_\alpha
    =
    e^{-\sigma_\alpha}\Omega,
\]
where the fourth equality is the local symplecticity of $F_{\alpha,h}$.
Hence $F_h^*\Omega=e^{\sigma_\alpha\circ F_h-\sigma_\alpha}\Omega$.
\end{proof}
\noindent
This is the discrete counterpart of the continuous relation
\eqref{eq:continuous-LCS-preservation}: the numerical map preserves the
defining LCS structure in the same conformal sense as the exact flow.

\paragraph{Independence of the local conformal potential.}

The exponent in \eqref{eq:global-discrete-LCS-preservation} is
independent of the choice of local potential on a connected component
of an overlap: if $\sigma_\alpha-\sigma_\beta=C_{\alpha\beta}$ on the
component containing both $z$ and $F_h(z)$, then
$\sigma_\alpha(F_h(z))-\sigma_\alpha(z)
=\sigma_\beta(F_h(z))-\sigma_\beta(z)$. Consequently,
\begin{equation}
    \rho_h(z)
    :=
    \exp
    \left(
        \sigma_\alpha(F_h(z))
        -
        \sigma_\alpha(z)
    \right)
    \label{eq:discrete-conformal-factor}
\end{equation}
defines a discrete conformal factor that is independent of the local
representation, provided that the initial point $z$ and the final point
$F_h(z)$ lie in the same connected component of the overlap of the
representations being compared, and the preservation law may be
written intrinsically as $F_h^*\Omega=\rho_h\,\Omega$. This proviso is
not vacuous: a discrete map may jump between two connected components
of $U_\alpha\cap U_\beta$ without leaving $U_\alpha$, in which case the
two local potentials differ by different constants at $z$ and at
$F_h(z)$. The one-step law then holds with the potential of the
representation in which the step is computed, which is the content of
Theorem~\ref{thm:global-LCS-preservation}.

\paragraph{Several time steps.}

The conformal preservation law has a multiplicative cocycle property
under iteration: for two steps,
$(F_h^2)^*\Omega=F_h^*(e^{\sigma_\alpha\circ F_h-\sigma_\alpha}\Omega)
=e^{\sigma_\alpha\circ F_h^2-\sigma_\alpha\circ F_h}
e^{\sigma_\alpha\circ F_h-\sigma_\alpha}\Omega$, and inductively
\begin{equation}
    (F_h^N)^*\Omega
    =
    e^{
        \sigma_\alpha\circ F_h^N
        -
        \sigma_\alpha
    }
    \Omega
    \label{eq:N-step-LCS-preservation}
\end{equation}
whenever the trajectory remains in a common representation.

More generally, when chart switching is required, composing the
one-step laws, each applied in the local representation in which the
corresponding step is computed, gives
\begin{equation}
    (F_h^N)^*\Omega
    =
    \exp
    \Big(
        \sum_{k=0}^{N-1}
        \big[
            \sigma_{\gamma_k}(z_{k+1})
            -
            \sigma_{\gamma_k}(z_k)
        \big]
    \Big)
    \Omega,
    \label{eq:N-step-LCS-preservation-switching}
\end{equation}
where $z_k=F_h^k(z)$ and $\gamma_k$ denotes the local representation
used for the $k$-th step. The exponent is a discrete line integral of
the Lee form along the numerical trajectory, and it is the discrete
counterpart of $\int_0^t\varphi(\xi_H)(\psi_s(x))\,\dd s$ in
\eqref{eq:continuous-conformal-factor-general}. If the representation
is changed from $\alpha$ to $\beta$ at $z_m$ and back to $\alpha$ at
$z_{m'}$, the exponent equals
\[
    \sigma_\alpha(z_N)
    -
    \sigma_\alpha(z_0)
    +
    C_{\alpha\beta}(z_m)
    -
    C_{\alpha\beta}(z_{m'}),
\]
where $C_{\alpha\beta}(z)$ denotes the value of the constant
$\sigma_\alpha-\sigma_\beta$ on the connected component of
$U_\alpha\cap U_\beta$ containing $z$. Thus
\eqref{eq:N-step-LCS-preservation} remains valid under chart switching
precisely when the switching points lie in the same connected component
of the overlap. When they do not, the additive constants do not cancel,
and the exponent in \eqref{eq:N-step-LCS-preservation-switching}
differs from $\sigma_\alpha(z_N)-\sigma_\alpha(z_0)$ by a period of
the Lee form. This is not a deficiency of the discretization but a
faithful reproduction of the continuous law
\eqref{eq:continuous-conformal-factor-general}, whose exponent along a
closed trajectory is $\oint\varphi$; the example with a non-exact Lee
form in Section~\ref{sec:numerics} exhibits this monodromy explicitly.

Two logically distinct properties therefore enter the global
preservation theorem: local structure preservation,
$F_{\alpha,h}^*\Omega_\alpha=\Omega_\alpha$, and overlap
compatibility,
$F_{\beta,h}\circ T_{\alpha\beta}=T_{\alpha\beta}\circ F_{\alpha,h}$.
Neither implies the other: two independently chosen symplectic
integrators may each preserve their respective local symplectic forms
while failing to describe the same discrete dynamics on an overlap,
and equivariance alone does not imply symplecticity. The numerical
experiments of Section~\ref{sec:numerics} measure the two properties
separately.

\paragraph{Order under globalization.}

Compatibility does not alter the local approximation order.

\begin{proposition}
\label{prop:order-globalization}
Suppose that, for each $\alpha$, the local discrete right Hamiltonian
$H_{\alpha,h}^{d+}$ approximates the corresponding exact discrete right
Hamiltonian $H_{\alpha,h}^{+,E}$ to order $p+1$ in the $C^1$ sense,
and suppose that the family $\{H_{\alpha,h}^{d+}\}$ is LCS-compatible.
Then the reconstructed global map $F_h$ has order $p$ with respect to
the exact LCS Hamiltonian flow.
\end{proposition}

\begin{proof}
By Theorem~\ref{thm:discrete-Hamiltonian-order}, each local map
$F_{\alpha,h}$ has order $p$ with respect to the exact local flow
$\Phi_{\alpha,h}$. Since $F_h=\Psi_\alpha^{-1}\circ F_{\alpha,h}\circ\Psi_\alpha$
and $\psi_h=\Psi_\alpha^{-1}\circ\Phi_{\alpha,h}\circ\Psi_\alpha$ with
$\Psi_\alpha$ smooth and independent of $h$, conjugation preserves the
local order, $F_h-\psi_h=\cO(h^{p+1})$ in local coordinates, and
compatibility ensures that these local representations define the same
global map on overlaps.
\end{proof}

\paragraph{A practical sufficient-condition principle.}

The preceding results may be summarized as a construction principle.

\begin{theorem}[Compatible Hamiltonian discretization principle]
\label{thm:compatible-discretization-principle}
Consider a cotangent-type LCS Hamiltonian system admitting local
canonical representations related on overlaps by
\[
    T_{\alpha\beta}(q,r)
    =
    (q,c_{\alpha\beta}r).
\]
Suppose that, on every chart, a regular discrete right Hamiltonian
$H_{\alpha,h}^{d+}$ is constructed such that, on every overlap,
\[
    H_{\beta,h}^{d+}(q,s)
    =
    c_{\alpha\beta}
    H_{\alpha,h}^{d+}
    \left(
        q,c_{\alpha\beta}^{-1}s
    \right).
\]
Then the local maps are equivariant with respect to the transition
maps, glue to a well-defined discrete dynamics $F_h$, and the resulting
global map satisfies
\begin{equation}
    F_h^*\Omega
    =
    e^{\sigma_\alpha\circ F_h-\sigma_\alpha}\Omega.
\end{equation}
If, in addition,
$H_{\alpha,h}^{d+}$ approximates $H_{\alpha,h}^{+,E}$ to order
$p+1$ in the $C^1$ sense, then $F_h$ is an order-$p$ approximation of
the exact LCS Hamiltonian flow.
\end{theorem}

\begin{proof}
The maps generated by the regular discrete right Hamiltonians
are locally symplectic by
Proposition~\ref{prop:local-discrete-symplecticity}. The conclusions then follow by combining
Theorems~\ref{thm:discrete-equivariance} and
\ref{thm:global-LCS-preservation},
Lemma~\ref{lem:discrete-gluing}, and
Proposition~\ref{prop:order-globalization}.
\end{proof}

\noindent
Theorem~\ref{thm:compatible-discretization-principle} reduces the
globalization problem to a concrete question at the level of the
numerical construction: how can one build approximate discrete
Hamiltonians that inherit the exact overlap transformation law? The
next section addresses this question for important classes of
geometric integrators.


\section{Classes of compatible structure-preserving integrators}
\label{sec:compatible-methods}

The results of Section~\ref{sec:globalization} reduce the construction
of global LCS integrators to two requirements. On each local canonical
representation, the numerical method must preserve the local symplectic
form, and the corresponding local maps must be equivariant with respect
to the overlap transformations $T_{\alpha\beta}$ of
\eqref{eq:T-alpha-beta}. In this section, we discuss classes of
numerical methods for which these requirements can be verified
directly. We first consider Hamiltonian variational constructions based
on approximations of the Type~II action. We then formulate a map-level
criterion that applies to symplectic Runge--Kutta methods: since the
overlap transformations are linear with constant coefficients,
Runge--Kutta equivariance follows from the corresponding equivariance
of the local Hamiltonian vector fields. Runge--Kutta-type methods and
their canonical and variational properties have been studied in
\cite{Suris89-RK, Suris90-RK}. This linearity (or, more generally,
affinity, see Section~\ref{subsec:affine}) is essential; the
statements below are not intended to assert equivariance of
Runge--Kutta methods under arbitrary nonlinear changes of local
coordinates.

\subsection{Compatible Type II Hamiltonian variational constructions}

A practical discrete Hamiltonian can be obtained by replacing the exact
trajectory and the exact action integral in
\eqref{eq:exact-discrete-right-Hamiltonian} by finite-dimensional
approximations. A general construction proceeds by choosing (i) a
finite-dimensional approximation space for the position curve; (ii) a
corresponding approximation of the momentum curve; (iii) a quadrature
rule for the Hamiltonian phase-space action; (iv) a discrete
stationarity condition determining the internal variables. Following
standard usage, $s$ denotes the number of stages; it should not be
confused with the momentum variable of the $\beta$-representation. Let
$0<\tau_1,\ldots,\tau_s\leq 1$ be quadrature nodes with weights
$b_1,\ldots,b_s$, and denote the approximate position and momentum
values at the quadrature nodes by $Q_i\approx q(\tau_i h)$ and
$R_i\approx r(\tau_i h)$. A generic approximation of the Type~II action
has the form
\begin{equation}
    \mathcal A_{\alpha,h}^{d}
    =
    r_1\cdot q_1
    -
    h
    \sum_{i=1}^{s}
    b_i
    \left[
        R_i\cdot V_i
        -
        H_\alpha(Q_i,R_i)
    \right],
    \label{eq:sec5-discrete-action-approx}
\end{equation}
where $V_i$ approximates $\dot q(\tau_i h)$. The associated discrete right
Hamiltonian is obtained by imposing the chosen finite-dimensional
stationarity conditions and eliminating the internal variables,
$H_{\alpha,h}^{d+}(q_0,r_1)=\operatorname{ext}\mathcal A_{\alpha,h}^{d}$,
where $\operatorname{ext}$ denotes evaluation at a stationary point of
the internal discrete variables, rather than necessarily a minimum. The
overlap scaling acts on momentum-type quantities by
$r\mapsto s=c_{\alpha\beta}r$, while the position variables are
unchanged, and the continuous Hamiltonian transforms according to
\eqref{eq:local-Hamiltonian-overlap}. Consequently, if the discrete
approximation space and the stationarity procedure are preserved by the
constant scaling of all momentum-type variables, the approximate action
inherits the same homogeneity.

\begin{proposition}
\label{prop:compatible-Galerkin}
Suppose that the discrete Type~II construction in each local
representation uses the same position approximation space, the same
quadrature rule, and momentum variables that transform according to
$S_i=c_{\alpha\beta}R_i$ and $s_1=c_{\alpha\beta}r_1$, and that the
continuous Hamiltonians satisfy \eqref{eq:local-Hamiltonian-overlap}.
Assume also that the stationary values defining the discrete
Hamiltonians are taken on corresponding regular branches related
by the same momentum scaling.
Then
\begin{equation}
    H_{\beta,h}^{d+}(q_0,s_1)
    =
    c_{\alpha\beta}
    H_{\alpha,h}^{d+}
    \left(
        q_0,c_{\alpha\beta}^{-1}s_1
    \right).
    \label{eq:sec5-Galerkin-compatibility}
\end{equation}
Hence the resulting discrete Hamiltonian maps are LCS-compatible.
\end{proposition}

\begin{proof}
Since the position variables and their discrete velocities are
unchanged, $S_i\cdot V_i=c_{\alpha\beta}R_i\cdot V_i$ and
$H_\beta(Q_i,S_i)=c_{\alpha\beta}H_\alpha(Q_i,R_i)$, so every term in
the discrete phase-space action scales by the same constant,
\begin{equation}
    \mathcal A_{\beta,h}^{d}
    =
    c_{\alpha\beta}
    \mathcal A_{\alpha,h}^{d}.
    \label{eq:sec5-action-scaling}
\end{equation}
Since \eqref{eq:sec5-action-scaling} holds identically in the internal
variables, $(Q^*,c_{\alpha\beta}R^*)$ is a stationary point of
$\mathcal A^d_{\beta,h}$ if and only if $(Q^*,R^*)$ is a stationary
point of $\mathcal A^d_{\alpha,h}$, so that stationary internal
variables correspond under the same momentum scaling. Evaluating the two
actions at corresponding stationary points gives
$H_{\beta,h}^{d+}(q_0,c_{\alpha\beta}r_1)
=c_{\alpha\beta}H_{\alpha,h}^{d+}(q_0,r_1)$, which is equivalent to
\eqref{eq:sec5-Galerkin-compatibility}.
\end{proof}
\noindent
By Theorem~\ref{thm:compatible-discretization-principle}, any regular
construction satisfying Proposition~\ref{prop:compatible-Galerkin}
produces a globally defined LCS structure-preserving integrator.

\subsection{Equivariance of Runge--Kutta discretizations}

A second important class of compatible methods is provided by
Runge--Kutta discretizations of the local Hamiltonian vector fields.
Let $\dot z=f_\alpha(z)$, $z=(q,r)$, denote the local Hamiltonian
system in the $\alpha$-representation and $\dot w=f_\beta(w)$,
$w=T_{\alpha\beta}z$, the corresponding system in the
$\beta$-representation. By
\eqref{eq:continuous-vector-field-equivariance},
\begin{equation}
    f_\beta(T_{\alpha\beta}z)
    =
    T_{\alpha\beta}f_\alpha(z),
    \label{eq:sec5-vector-field-equivariance}
\end{equation}
where, since $T_{\alpha\beta}$ is linear, the same notation denotes the
map and its tangent map. Consider an $s$-stage Runge--Kutta method
with coefficients $(a_{ij},b_i)$, whose stages and update in the
$\alpha$-representation are
\begin{equation}
    Z_i
    =
    z_n
    +
    h
    \sum_{j=1}^{s}
    a_{ij}f_\alpha(Z_j),
    \qquad
    z_{n+1}
    =
    z_n
    +
    h
    \sum_{i=1}^{s}
    b_i f_\alpha(Z_i),
    \label{eq:sec5-RK-stages-alpha}
\end{equation}
and similarly in the $\beta$-representation with $f_\beta$.

\begin{proposition}[Runge--Kutta equivariance]
\label{prop:RK-equivariance}
Suppose that the local vector fields satisfy
$f_\beta\circ T_{\alpha\beta}=T_{\alpha\beta}\circ f_\alpha$, where
$T_{\alpha\beta}$ is linear, and suppose that the same Runge--Kutta
tableau is used in both local representations. Whenever the stage
equations determine corresponding local solutions, the Runge--Kutta
maps satisfy
\begin{equation}
    F_{\beta,h}^{RK}\circ T_{\alpha\beta}
    =
    T_{\alpha\beta}\circ F_{\alpha,h}^{RK}.
    \label{eq:sec5-RK-equivariance}
\end{equation}
\end{proposition}

\begin{proof}
Applying the linear map $T_{\alpha\beta}$ to the stage and update
equations \eqref{eq:sec5-RK-stages-alpha} and using
$T_{\alpha\beta}f_\alpha(Z_j)=f_\beta(T_{\alpha\beta}Z_j)$ shows that
$W_i=T_{\alpha\beta}Z_i$ and $w_{n+1}=T_{\alpha\beta}z_{n+1}$ satisfy
the stage and update equations of the same method for $f_\beta$ with
initial value $w_n=T_{\alpha\beta}z_n$; by uniqueness of the stage
solution,
$F_{\beta,h}^{RK}(T_{\alpha\beta}z_n)=T_{\alpha\beta}(F_{\alpha,h}^{RK}(z_n))$.
\end{proof}

\begin{remark}
The linearity assumption in
Proposition~\ref{prop:RK-equivariance} cannot be dropped in favour of
general nonlinear transitions. Standard Runge--Kutta methods are affine
equivariant, but they are not generally equivariant under arbitrary
nonlinear coordinate transformations. The cotangent-type overlap maps
considered here are constant linear fiber scalings, and hence lie
within the class for which the preceding argument applies directly;
the affine case is treated in Proposition~\ref{prop:RK-affine-equivariance}.
\end{remark}

\paragraph{Symplectic Runge--Kutta methods.}

Equivariance provides the overlap compatibility required for
globalization. To obtain an LCS structure-preserving method, the local
Runge--Kutta discretization must additionally be symplectic. For a
canonical Hamiltonian system, an $s$-stage Runge--Kutta method is
symplectic if its coefficients satisfy
\begin{equation}
    b_i a_{ij}
    +
    b_j a_{ji}
    =
    b_i b_j,
    \qquad
    i,j=1,\ldots,s,
    \label{eq:sec5-symplectic-RK-condition}
\end{equation}
in which case $(F_{\alpha,h}^{RK})^*\Omega_\alpha=\Omega_\alpha$
\cite{HaLuWa}. Combining this fact with
Proposition~\ref{prop:RK-equivariance} immediately gives the following
consequence.

\begin{theorem}[Global LCS integrators from symplectic Runge--Kutta methods]
\label{thm:symplectic-RK-LCS}
Consider the cotangent-type LCS setting of
Section~\ref{sec:continuous}, with overlap transformations
$T_{\alpha\beta}(q,r)=(q,c_{\alpha\beta}r)$, and apply the same
symplectic Runge--Kutta method to every local canonical Hamiltonian
system. Then the resulting local maps are equivariant on overlaps and
glue to a global discrete map $F_h$ satisfying
\begin{equation}
    F_h^*\Omega
    =
    e^{\sigma_\alpha\circ F_h-\sigma_\alpha}\Omega.
    \label{eq:sec5-RK-global-LCS}
\end{equation}
\end{theorem}

\begin{proof}
The local Hamiltonian vector fields are equivariant with respect to the
linear maps $T_{\alpha\beta}$ by
\eqref{eq:sec5-vector-field-equivariance}, so
Proposition~\ref{prop:RK-equivariance} gives equivariance of the local
Runge--Kutta maps, and \eqref{eq:sec5-symplectic-RK-condition} gives
local symplecticity. The conclusion follows from
Lemma~\ref{lem:discrete-gluing} and
Theorem~\ref{thm:global-LCS-preservation}.
\end{proof}

\subsection{Collocation and Gauss--Legendre methods}
\label{subsec:gauss-legendre}

The family of methods that realizes
Theorem~\ref{thm:symplectic-RK-LCS} in its sharpest form is the
Gauss--Legendre collocation family \cite{HaLuWa,SaCa}. Let $s\geq1$
and let $0<\tau_1<\cdots<\tau_s<1$ be the roots of the shifted Legendre
polynomial
$\widetilde P_s(\lambda)=\frac{1}{s!}\frac{\dd^s}{\dd\lambda^s}[\lambda^s(\lambda-1)^s]$
on $[0,1]$. With $\ell_1,\ldots,\ell_s$ the Lagrange interpolation
polynomials associated with these nodes, the coefficients
\begin{equation}
    a_{ij}
    =
    \int_0^{\tau_i}\ell_j(\lambda)\,\dd\lambda,
    \qquad
    b_i
    =
    \int_0^{1}\ell_i(\lambda)\,\dd\lambda,
    \label{eq:GL-coefficients}
\end{equation}
define the \emph{$s$-stage Gauss--Legendre method}. Equivalently, it is
the collocation method that seeks a polynomial $u$ of degree $s$ with
$u(0)=z_n$ and $\dot u(\tau_ih)=f_\alpha(u(\tau_ih))$, $i=1,\ldots,s$, and
sets $z_{n+1}=u(h)$. The first three members are
\[
\begin{array}{c|c}
\frac12 & \frac12\\
\hline
& 1
\end{array}
\qquad
\begin{array}{c|cc}
\frac12-\frac{\sqrt3}{6}
& \frac14
& \frac14-\frac{\sqrt3}{6}
\\
\frac12+\frac{\sqrt3}{6}
& \frac14+\frac{\sqrt3}{6}
& \frac14
\\
\hline
& \frac12 & \frac12
\end{array}
\qquad
\begin{array}{c|ccc}
\frac12-\frac{\sqrt{15}}{10}
& \frac{5}{36}
& \frac{2}{9}-\frac{\sqrt{15}}{15}
& \frac{5}{36}-\frac{\sqrt{15}}{30}
\\
\frac12
& \frac{5}{36}+\frac{\sqrt{15}}{24}
& \frac{2}{9}
& \frac{5}{36}-\frac{\sqrt{15}}{24}
\\
\frac12+\frac{\sqrt{15}}{10}
& \frac{5}{36}+\frac{\sqrt{15}}{30}
& \frac{2}{9}+\frac{\sqrt{15}}{15}
& \frac{5}{36}
\\
\hline
& \frac{5}{18} & \frac{4}{9} & \frac{5}{18}
\end{array}
\]
corresponding to the implicit midpoint rule ($s=1$) and to the
fourth- and sixth-order Gauss methods.
We collect the classical facts that will be used below.

\begin{proposition}
\label{prop:GL-classical}
The $s$-stage Gauss--Legendre method has order $2s$, its underlying
quadrature rule being exact for polynomials of degree $\leq2s-1$ and
the collocation polynomial superconverging at the endpoint. It is
symmetric, $A$-stable and algebraically stable, with stability function
the $(s,s)$ diagonal Pad\'e approximation of the exponential, and it
preserves exactly all quadratic first integrals of the differential
equation. Its coefficients satisfy
\begin{equation}
    b_ia_{ij}+b_ja_{ji}-b_ib_j=0,
    \qquad
    i,j=1,\ldots,s,
    \label{eq:GL-symplectic-condition}
\end{equation}
so that \eqref{eq:sec5-symplectic-RK-condition} holds and the method is
symplectic for canonical Hamiltonian systems.
Finally, if $f_\alpha$ is smooth, then for each initial point
$z_*$ the stage equations admit a unique local solution branch
near $Z_1=\cdots=Z_s=z_*$ for $(z_n,h)$ sufficiently close
to $(z_*,0)$. This branch depends smoothly on $z_n$ and $h$.
To see this, write the stage equations as
\[
    \mathcal G_i(Z,z_n,h)
    =
    Z_i-z_n-h\sum_{j=1}^{s}a_{ij}f_\alpha(Z_j)=0,
\]
where $Z=(Z_1,\ldots,Z_s)$ and
$\mathcal G=(\mathcal G_1,\ldots,\mathcal G_s)$. At $h=0$, these equations have the solution $Z_i=z_n$
for all $i$, and the Jacobian $D_Z\mathcal G$ at this solution
is the identity. The claim then follows from the implicit
function theorem.
\end{proposition}
\noindent
This local solvability result ensures that the one-step maps
$F_{\alpha,h}$ are well defined and smooth on the selected
stage branches for sufficiently small $h$.
By Proposition~\ref{prop:RK-equivariance}, stage solutions
are carried bijectively to stage solutions in corresponding
overlap neighborhoods. Consequently, local existence and
uniqueness transfer between representations at the same
step size, provided the stages remain in these neighborhoods.
The same conclusion holds for affine transitions by
Proposition~\ref{prop:RK-affine-equivariance}.

\subsubsection*{Gauss--Legendre as a Type II Galerkin construction}

Gauss--Legendre methods occupy a distinguished position in the present
framework because they belong simultaneously to both classes of
compatible methods identified in this section: they arise from a
Type~II Galerkin approximation of the Hamiltonian phase-space action,
and they are symplectic Runge--Kutta methods.

Fix a local representation and consider the Type~II action
\eqref{eq:sec5-discrete-action-approx}. Choose
\begin{enumerate}[label=(\roman*)]
    \item as position approximation space the polynomials
    $q_d$ of degree $\leq s$ on $[0,h]$ with $q_d(0)=q_0$;
    \item as momentum approximation the values $R_i$ at the quadrature
    nodes, i.e., an unconstrained polynomial $r_d$ of degree $\leq s-1$,
    the final momentum $r_1$ entering only through the boundary term;
    \item the $s$-point Gauss--Legendre quadrature rule with nodes and
    weights \eqref{eq:GL-coefficients}.
\end{enumerate}
\noindent
Writing $Q_i=q_d(\tau_ih)$, $R_i=r_d(\tau_ih)$ and $V_i=\dot q_d(\tau_ih)$, the
discrete right Hamiltonian is
\begin{equation}
    H_{\alpha,h}^{d+}(q_0,r_1)
    =
    \operatorname{ext}
    \left\{
        r_1\cdot q_1
        -
        h\sum_{i=1}^{s}b_i
        \left[
            R_i\cdot V_i-H_\alpha(Q_i,R_i)
        \right]
    \right\},
    \label{eq:GL-discrete-right-H}
\end{equation}
the extremum being taken over the internal variables of $q_d$, namely
its nodal velocities $V_j$, and over the nodal momenta $R_i$, with
$q_1=q_d(h)$.

\begin{proposition}
\label{prop:GL-is-variational}
The discrete right Hamiltonian \eqref{eq:GL-discrete-right-H}
generates, through the right discrete Hamilton equations
\eqref{eq:right-discrete-Hamilton-equations}, the $s$-stage
Gauss--Legendre method applied to the local canonical Hamiltonian
system associated with $H_\alpha$. Consequently, the induced local map
has order $2s$, and \eqref{eq:GL-discrete-right-H} approximates the
exact discrete right Hamiltonian $H_{\alpha,h}^{+,E}$ to order
$2s+1$ in the $C^1$ sense, up to an additive constant.
\end{proposition}

\begin{proof}
We spell out the correspondence because it is the point at which the
Type~II variational construction and the Runge--Kutta description meet.
Let
\[
    Q_i=q_d(\tau_i h),
    \qquad
    R_i=r_d(\tau_i h),
    \qquad
    i=1,\ldots,s,
\]
and write the canonical Hamiltonian vector field as
\[
    \dot q=D_2H_\alpha(q,r),
    \qquad
    \dot r=-D_1H_\alpha(q,r).
\]
We use the nodal velocities $V_j=\dot q_d(\tau_jh)$,
$j=1,\ldots,s$, as internal variables parametrizing $q_d$
with fixed initial position $q_0$. Stationarity of
\eqref{eq:GL-discrete-right-H} with respect to the
nodal momenta $R_i$ gives $V_i=D_2H_\alpha(Q_i,R_i)$, and stationarity
with respect to the nodal velocities $V_j$ gives
$R_j=r_1+h\sum_i(b_ia_{ij}/b_j)\,D_1H_\alpha(Q_i,R_i)$; together with
$r_0=D_1H^{d+}_{\alpha,h}(q_0,r_1)=r_1+h\sum_ib_iD_1H_\alpha(Q_i,R_i)$,
this gives
\[
    R_j
    =
    r_0
    -
    h\sum_{i=1}^{s}
    \Bigl(b_i-\frac{b_ia_{ij}}{b_j}\Bigr)
    D_1H_\alpha(Q_i,R_i),
\]
and these relations show, as in \cite{LeZh11}, that the stage values satisfy
a symplectic partitioned Runge--Kutta system of the form
\begin{align}
    Q_i
    &=
    q_0
    +h\sum_{j=1}^{s}a_{ij}
    D_2H_\alpha(Q_j,R_j),
    \label{eq:GL-stage-Q}
    \\
    R_i
    &=
    r_0
    -h\sum_{j=1}^{s}\widetilde a_{ij}
    D_1H_\alpha(Q_j,R_j),
    \label{eq:GL-stage-R}
\end{align}
where the coefficients conjugate to $a_{ij}$ are
\begin{equation}
    \widetilde a_{ij}
    =
    b_j-\frac{b_j}{b_i}a_{ji}.
    \label{eq:GL-conjugate-coefficients}
\end{equation}
The endpoint values are
\begin{align}
    q_1
    &=
    q_0
    +h\sum_{i=1}^{s}b_iD_2H_\alpha(Q_i,R_i),
    \label{eq:GL-endpoint-q}
    \\
    r_1
    &=
    r_0
    -h\sum_{i=1}^{s}b_iD_1H_\alpha(Q_i,R_i).
    \label{eq:GL-endpoint-r}
\end{align}
For the Gauss--Legendre coefficients,
\eqref{eq:GL-symplectic-condition} gives
\[
    b_i a_{ij}+b_j a_{ji}=b_i b_j,
\]
and therefore
\[
    \widetilde a_{ij}
    =
    b_j-\frac{b_j}{b_i}a_{ji}
    =a_{ij}.
\]
Thus the partitioned scheme collapses to the ordinary $s$-stage
Gauss--Legendre Runge--Kutta method applied to the canonical Hamiltonian
system. In particular, the map obtained from the Type~II generating
function is the same map as the Gauss collocation map considered in
Proposition~\ref{prop:GL-classical}.

For the order statement, the induced map is the $s$-stage Gauss
collocation map, which has order $2s$ by
Proposition~\ref{prop:GL-classical}. By
Remark~\ref{rem:order-converse}, \eqref{eq:GL-discrete-right-H}
therefore approximates $H^{+,E}_{\alpha,h}$ to order $2s+1$ in the
$C^1$ sense, up to an additive constant, and
Theorem~\ref{thm:discrete-Hamiltonian-order} with $p=2s$ returns the
order $2s$ of the induced map.
\end{proof}

\begin{remark}
Proposition~\ref{prop:GL-is-variational} is the reason why the two
routes developed above---the generating-function route and the
Runge--Kutta equivariance route---are not merely parallel. For
Gauss--Legendre methods the verification of the
compatibility condition \eqref{eq:compatible-discrete-H} can be carried
out either at the level of the generating function, through
Proposition~\ref{prop:compatible-Galerkin}, or at the level of the
vector field, through Proposition~\ref{prop:RK-equivariance}, and the
two verifications concern the same numerical method.
\end{remark}

\subsubsection*{Compatibility of the Gauss--Legendre family}

\begin{proposition}
\label{prop:GL-compatible}
Let the same $s$-stage Gauss--Legendre tableau be used in every local
canonical representation. Then the family of generated discrete right
Hamiltonians is LCS-compatible in the sense of
Definition~\ref{def:LCS-compatible-discrete-H}, and the local stage
variables transform as
\begin{equation}
    Q_i^\beta=Q_i^\alpha,
    \qquad
    S_i=c_{\alpha\beta}R_i.
    \label{eq:GL-stage-scaling}
\end{equation}
\end{proposition}

\begin{proof}
The Gauss quadrature and the polynomial approximation spaces are
identical in both representations, the position variables are
unaffected by $T_{\alpha\beta}$ and all momentum-type variables scale
by $c_{\alpha\beta}$, so Proposition~\ref{prop:compatible-Galerkin}
applies and yields \eqref{eq:compatible-discrete-H}. The stage scaling
\eqref{eq:GL-stage-scaling} is Proposition~\ref{prop:RK-equivariance}
applied to the stages \eqref{eq:GL-stage-Q}--\eqref{eq:GL-stage-R},
which by Proposition~\ref{prop:GL-is-variational} are those of the
Gauss--Legendre method: setting $s_n=c_{\alpha\beta}r_n$ and
$S_i=c_{\alpha\beta}R_i$ and using \eqref{eq:H-overlap-D1}, the
position stage equations are unchanged and the momentum stage equations
are multiplied by $c_{\alpha\beta}$.
\end{proof}

\begin{corollary}
\label{cor:GL-LCS-integrator}
Under the hypotheses of
Theorem~\ref{thm:compatible-discretization-principle}, the $s$-stage
Gauss--Legendre method applied in every local canonical representation
defines a global discrete map $F_h$ on the LCS manifold satisfying
\begin{equation}
    F_h^*\Omega
    =
    e^{\sigma_\alpha\circ F_h-\sigma_\alpha}\Omega,
\end{equation}
and $F_h$ is an approximation of order $2s$ of the LCS Hamiltonian
flow. In particular the framework contains a hierarchy of
arbitrarily high even-order LCS structure-preserving integrators.
\end{corollary}
The following remark explains why, on compact sets and for sufficiently
small step sizes, each numerical step and its stages can be computed
within a single chart.
\begin{remark}[Steps near chart boundaries]
\label{rem:chart-boundaries}
Let $K\subset M$ be compact and let $d$ be a Riemannian distance.
Choose a finite cover of $K$ by open sets whose closures are
compact and contained in connected chart domains. There exists
$\delta>0$ such that, for each $z\in K$, a corresponding chart
domain $U_\alpha$ contains $B_d(z,\delta)$.

For the consistent local solution branch, the endpoint and all
Runge--Kutta stage points $Z_i$, $i=1,\ldots,s$, equal $z$ at
$h=0$. Smooth dependence on $(z,h)$ and compactness therefore
give a common constant $C>0$ such that
\[
    d(z,F_{\alpha,h}(z))\leq Ch,
    \qquad
    \max_{1\leq i\leq s}d(z,Z_i)\leq Ch
\]
for sufficiently small $h$. Choosing $Ch<\delta$ ensures that
the initial point, the endpoint, and all stages lie in
$U_\alpha$.
Since $B_d(z,\delta)$ is path connected, whenever it is contained in
two chart domains it lies in a single connected component of their
overlap, so a single constant $c_{\alpha\beta}$ applies and overlap
compatibility gives the same update.
\end{remark}

\subsubsection*{Modified Hamiltonians and long-time behaviour}

A further consequence of symplecticity is of direct relevance to the
numerical experiments. In each local representation, backward error
analysis \cite{HaLuWa,Reich-99} provides a formal modified Hamiltonian
\begin{equation}
    \widetilde H_{\alpha,h}
    =
    H_\alpha
    +
    h^{2s}H_\alpha^{[2s]}
    +
    h^{2s+2}H_\alpha^{[2s+2]}
    +\cdots.
    \label{eq:modified-Hamiltonian}
\end{equation}
For analytic $H_\alpha$ and sufficiently small $h$, assuming that the
numerical orbit remains in a fixed compact subset of the chart with
uniform analytic control and on a common regular stage branch,
a modified Hamiltonian optimally truncated at an order
proportional to $h^{-1}$ has an exponentially small
one-step conservation defect \cite[Ch.~IX]{HaLuWa}.
Consequently, $H_\alpha$ is conserved up to
$\mathcal O(h^{2s})$ over exponentially long times,
as long as these assumptions remain valid.

\vspace{3mm}

\noindent
Two remarks are relevant for the LCS setting.

\begin{remark}
\label{rem:H-conformal-evolution}
The function $H$ is not, in general, a first integral of the
LCS dynamics. Contracting
\eqref{eq:LCS-Hamiltonian-vector-field} with $\xi_H$ gives
\[
    \xi_H(H)=H\varphi(\xi_H),
\]
so that, along an exact trajectory remaining in $U_\alpha$,
\begin{equation}
    H\circ\psi_t
    =
    e^{\sigma_\alpha\circ\psi_t-\sigma_\alpha}H.
    \label{eq:H-conformal-law}
\end{equation}
Thus, $H$ obeys the same conformal scaling law as $\Omega$,
and $H_\alpha=e^{-\sigma_\alpha}H$ is conserved. For the $s$-stage
Gauss--Legendre method, under the backward error analysis assumptions
stated above, a numerical trajectory $z_n=F_h^n(z_0)$ remaining in a
fixed compact subset of $U_\alpha$ satisfies
\begin{equation}
    e^{-\sigma_\alpha(z_n)}H(z_n)
    =
    e^{-\sigma_\alpha(z_0)}H(z_0)
    +
    \mathcal O(h^{2s}),
    \label{eq:discrete-H-conformal-law}
\end{equation}
uniformly in $n$ over exponentially long time intervals. Equivalently,
whenever $H_\alpha(z_0)\neq0$,
\[
    H(z_n)
    =
    e^{\sigma_\alpha(z_n)-\sigma_\alpha(z_0)}H(z_0)
    \left(1+\mathcal O(h^{2s})\right).
\]
The uniform bound therefore holds for the local Hamiltonian, or in
relative form for $H$, but not for the absolute error in $H$ itself,
which carries the factor $e^{\sigma_\alpha(z_n)}$ and may grow along
the trajectory. Under chart switching, the accumulated conformal factor
is given by \eqref{eq:N-step-LCS-preservation-switching};
for closed paths, it records the periods of the Lee form
and hence its monodromy.
The long-time energy estimate above is restricted to
trajectories remaining in a fixed compact subset of
a single chart.
The weighted relation \eqref{eq:discrete-H-conformal-law}
is the local energy-type prediction tested numerically.
\end{remark}

\begin{remark}
\label{rem:modified-equivariance}
The modified Hamiltonians are themselves compatible. In each local
representation, the term $H_\alpha^{[k]}$ in
\eqref{eq:modified-Hamiltonian} is a universal linear combination,
with coefficients determined by the method alone and with no additive
constants, of elementary Hamiltonians built from $k+1$ factors given by
$H_\alpha$ and its derivatives, contracted by means of $k$ copies of
the Poisson tensor $\Omega_\alpha^{-1}$. Under the transition map one
has $H_\beta=c_{\alpha\beta}H_\alpha\circ T_{\alpha\beta}^{-1}$, while,
by $T_{\alpha\beta}^*\Omega_\beta
    =
    c_{\alpha\beta}\Omega_\alpha$, the Poisson tensor of
$\Omega_\beta$ is the push-forward by $T_{\alpha\beta}$ of
$c_{\alpha\beta}^{-1}\Omega_\alpha^{-1}$. Each elementary Hamiltonian
entering $H^{[k]}_\beta$ is therefore
$c_{\alpha\beta}^{k+1}c_{\alpha\beta}^{-k}=c_{\alpha\beta}$ times the
corresponding term of $H^{[k]}_\alpha$ composed with
$T_{\alpha\beta}^{-1}$, and one obtains
\[
    \widetilde H_{\beta,h}(q,s)
    =
    c_{\alpha\beta}
    \widetilde H_{\alpha,h}
    \left(q,c_{\alpha\beta}^{-1}s\right).
\]
Consequently the local functions $e^{\sigma_\alpha}\widetilde
H_{\alpha,h}$ agree on overlaps and define a global function
$\widetilde H_h$ on $M$. Hence the modified equations also glue, and
the numerical trajectory is exponentially close to the exact flow of
the nearby \emph{globally defined} LCS Hamiltonian system associated
with $\widetilde H_h$. This is the structural reason behind the
absence of drift observed in Section~\ref{sec:numerics}.
\end{remark}

\paragraph{The implicit midpoint rule.}

The simplest Gauss--Legendre method is the one-stage implicit midpoint
rule, a second-order symplectic Runge--Kutta method. For the local
canonical Hamiltonian system it reads
\begin{equation}
    \frac{q_{n+1}-q_n}{h}
    =
    D_2H_\alpha
    \left(
        \frac{q_n+q_{n+1}}{2},
        \frac{r_n+r_{n+1}}{2}
    \right),
    \qquad
    \frac{r_{n+1}-r_n}{h}
    =
    -
    D_1H_\alpha
    \left(
        \frac{q_n+q_{n+1}}{2},
        \frac{r_n+r_{n+1}}{2}
    \right),
    \label{eq:sec5-midpoint-q}
\end{equation}
and its equivariance under $T_{\alpha\beta}$ and the conformal law of
its reconstruction are the case $s=1$ of
Proposition~\ref{prop:GL-compatible} and
Corollary~\ref{cor:GL-LCS-integrator}; the general globalization
principle thus recovers the implicit midpoint construction as a
particular instance rather than as an isolated example. For the
running Example~\ref{ex:continuous-LCS-oscillator}, with midpoint
quantities $\bar q_n=\tfrac12(q_n+q_{n+1})$ and similarly for
$\bar r_n$, $\bar x_n$ and $\bar s_n$, the scheme is
\begin{equation}
\begin{aligned}
    q_{n+1}-q_n
    &=
    h e^{-\bar x_n}\bar r_n,
    \\[1mm]
    r_{n+1}-r_n
    &=
    -hk e^{\bar x_n}\bar q_n,
    \\[1mm]
    x_{n+1}-x_n
    &=
    h e^{-\bar x_n}\bar s_n,
    \\[1mm]
    s_{n+1}-s_n
    &=
    h
    \left\{
        \frac12e^{-\bar x_n}
        \left(
            \bar r_n^2+\bar s_n^2
        \right)
        -
        e^{\bar x_n}
        \left[
            \frac{k}{2}\bar q_n^2
            +
            \mu(\bar x_n+1)
        \right]
    \right\}.
\end{aligned}
\label{eq:sec5-example-midpoint}
\end{equation}
It preserves the local canonical symplectic form
$\Omega_{\mathrm{loc}}=\dd q\wedge\dd r+\dd x\wedge\dd s$, and since
$\sigma=-x$ and $\Omega=e^{-x}\Omega_{\mathrm{loc}}$, the reconstructed
one-step map satisfies
$(F_h^{MP})^*\Omega=e^{-(x_{n+1}-x_n)}\Omega$: the discrete conformal
multiplier is determined exactly by the change of the local conformal
potential across one numerical time step.

\subsection{Extension to affine base transitions}
\label{subsec:affine}

The preceding results were formulated for overlap maps
that preserve the base coordinates and rescale the momenta.
We now extend the relevant transformation and compatibility
results to affine base transitions. This extension includes
the angular chart transitions used in Example~\ref{ex:S1-example}.

\begin{proposition}[Exact discrete right Hamiltonians under affine transitions]
\label{prop:exact-discrete-H-affine}
Consider an affine cotangent transition
\begin{equation}
    T_{\alpha\beta}(q,r)
    =
    \bigl(Bq+b,c_{\alpha\beta}B^{-T}r\bigr),
    \label{eq:affine-transition}
\end{equation}
where $B$ is a constant invertible matrix, $b$ is a constant
vector, and $c_{\alpha\beta}>0$ is constant on the connected
overlap component under consideration. Suppose that the local
Hamiltonians satisfy
\[
    H_\beta\circ T_{\alpha\beta}
    =
    c_{\alpha\beta}H_\alpha.
\]
For mixed data $(q_0,r_1)$ such that the corresponding solution
of the $\alpha$-Hamiltonian system remains in this component
throughout $[0,h]$, set
\[
    Q_0=Bq_0+b,
    \qquad
    S_1=c_{\alpha\beta}B^{-T}r_1.
\]
On corresponding solution branches, the exact discrete right
Hamiltonians satisfy
\begin{equation}
    H_{\beta,h}^{+,E}(Q_0,S_1)
    =
    c_{\alpha\beta}H_{\alpha,h}^{+,E}
    \left(
        B^{-1}(Q_0-b),
        c_{\alpha\beta}^{-1}B^TS_1
    \right)
    +
    S_1\cdot b.
    \label{eq:exact-discrete-H-affine}
\end{equation}
\end{proposition}

\begin{proof}
The identities
$T_{\alpha\beta}^*\Omega_\beta=c_{\alpha\beta}\Omega_\alpha$
and
$H_\beta\circ T_{\alpha\beta}=c_{\alpha\beta}H_\alpha$
imply that $T_{\alpha\beta}$ maps $\alpha$-Hamiltonian
trajectories to $\beta$-Hamiltonian trajectories.
Thus, a solution $(q(t),r(t))$ with
$q(0)=q_0$ and $r(h)=r_1$ is transformed into
\[
    Q(t)=Bq(t)+b,
    \qquad
    S(t)=c_{\alpha\beta}B^{-T}r(t),
\]
with $Q(0)=Q_0$ and $S(h)=S_1$.
Along these corresponding trajectories,
\[
    S\cdot\dot Q-H_\beta(Q,S)
    =
    c_{\alpha\beta}
    \bigl(r\cdot\dot q-H_\alpha(q,r)\bigr),
\]
and the endpoint term satisfies
\[
    S_1\cdot Q(h)
    =
    c_{\alpha\beta}r_1\cdot q(h)+S_1\cdot b.
\]
Hence, by the definition of the exact discrete right Hamiltonian,
\[
\begin{aligned}
    H_{\beta,h}^{+,E}(Q_0,S_1)
    &=
    S_1\cdot Q(h)
    -
    \int_0^h
    \bigl(S\cdot\dot Q-H_\beta(Q,S)\bigr)\,\dd t\\
    &=
    c_{\alpha\beta}H_{\alpha,h}^{+,E}(q_0,r_1)
    +S_1\cdot b.
\end{aligned}
\]
Substituting $q_0=B^{-1}(Q_0-b)$ and
$r_1=c_{\alpha\beta}^{-1}B^TS_1$
gives \eqref{eq:exact-discrete-H-affine}.
\end{proof}

\begin{corollary}[Affine compatibility]
\label{cor:affine-TypeII-compatibility}
Suppose that regular discrete right Hamiltonians satisfy
\[
    H_{\beta,h}^{d+}(Q_0,S_1)
    =
    c_{\alpha\beta}H_{\alpha,h}^{d+}
    \left(
        B^{-1}(Q_0-b),
        c_{\alpha\beta}^{-1}B^TS_1
    \right)
    +S_1\cdot b.
\]
Then their generated maps satisfy
\[
    F_{\beta,h}\circ T_{\alpha\beta}
    =
    T_{\alpha\beta}\circ F_{\alpha,h}
\]
on corresponding local solution branches where both sides
are defined.
\end{corollary}

\begin{proof}
Differentiating the transformation law gives
\[
    D_1H_{\beta,h}^{d+}
    =
    c_{\alpha\beta}B^{-T}D_1H_{\alpha,h}^{d+},
    \qquad
    D_2H_{\beta,h}^{d+}
    =
    BD_2H_{\alpha,h}^{d+}+b,
\]
where the derivatives on the right are evaluated at
$(B^{-1}(Q_0-b),c_{\alpha\beta}^{-1}B^TS_1)$.
The discrete Hamilton equations therefore yield
$S_0=c_{\alpha\beta}B^{-T}r_0$ and $Q_1=Bq_1+b$,
which proves the claim.
\end{proof}
The affine equivariance of Runge--Kutta methods is a classical result; see, \cite{McLachlan-B-series}. For completeness, we give a direct proof in the present setting.
\begin{proposition}[Runge--Kutta equivariance under affine transitions]
\label{prop:RK-affine-equivariance}
Let
$T_{\alpha\beta}(z)=A_{\alpha\beta}z+d_{\alpha\beta}$,
where $A_{\alpha\beta}$ is a constant invertible matrix and
$d_{\alpha\beta}$ is a constant vector. Suppose that
\[
    f_\beta(T_{\alpha\beta}z)
    =
    A_{\alpha\beta}f_\alpha(z),
\]
and that the same Runge--Kutta tableau is used in both
representations. Whenever the stage equations admit unique
local solutions on corresponding branches, with all stages
lying in the relevant overlap domains, the Runge--Kutta maps
satisfy
\[
    F_{\beta,h}^{RK}\circ T_{\alpha\beta}
    =
    T_{\alpha\beta}\circ F_{\alpha,h}^{RK}.
\]
\end{proposition}

\begin{proof}
Let $Z_i$ be the stages in the $\alpha$-representation,
and set
$w_n=T_{\alpha\beta}z_n$ and
$W_i=T_{\alpha\beta}Z_i$.
Applying $A_{\alpha\beta}$ to the stage equations and
adding $d_{\alpha\beta}$ gives
\[
    W_i
    =
    A_{\alpha\beta}z_n+d_{\alpha\beta}
    +
    h\sum_{j=1}^{s}
    a_{ij}A_{\alpha\beta}f_\alpha(Z_j)
    =
    w_n+
    h\sum_{j=1}^{s}a_{ij}f_\beta(W_j).
\]
Thus, the $W_i$ satisfy the stage equations in the
$\beta$-representation. Similarly, the update satisfies
\[
    T_{\alpha\beta}z_{n+1}
    =
    A_{\alpha\beta}z_n+d_{\alpha\beta}
    +
    h\sum_{i=1}^{s}
    b_iA_{\alpha\beta}f_\alpha(Z_i)
    =
    w_n+
    h\sum_{i=1}^{s}b_i f_\beta(W_i)
    =
    w_{n+1}.
\]
Uniqueness on the corresponding stage branches therefore
yields the claimed equivariance.
\end{proof}

\begin{remark}[Order under affine transitions]
Affine transitions preserve the local error order.
Indeed, equivariance of the numerical and exact flows gives
\[
    F_{\beta,h}(T_{\alpha\beta}z)
    -
    \Phi_{\beta,h}(T_{\alpha\beta}z)
    =
    A_{\alpha\beta}
    \bigl(F_{\alpha,h}(z)-\Phi_{\alpha,h}(z)\bigr).
\]
Moreover, if the approximate and exact discrete right
Hamiltonians satisfy the same affine transformation law,
their endpoint terms cancel upon subtraction.
The remaining difference is obtained by constant scaling
and affine composition, which preserve its $C^1$
approximation order on corresponding compact sets.
Thus, the preceding local order results apply under
the same regularity assumptions.
\end{remark}

The Hamiltonian construction does not require passage through a
Lagrangian formulation. Nevertheless, when a regular Lagrangian
description is available, the same overlap structure is inherited by
the continuous and discrete Legendre transformations. We develop this
counterpart in the next section.


\section{The Lagrangian counterpart}
\label{sec:lagrangian}

The globalization results of Section~\ref{sec:globalization} were
obtained on the Hamiltonian side and therefore did not require a
Legendre transformation. When a regular Lagrangian description is
available, the same local-to-global structure admits a variational
formulation in terms of discrete Lagrangians; see
\cite{guibout2004discrete,LeeLeokMc,MaDiMa06,MarsWest, ScShLe-2018,Wend-Marsden}
for variational integrators for Lagrangian systems and
\cite{Colombo-2023, deLeon-time-discrete} for
time-dependent Lagrangian systems. We state the continuous and discrete
Lagrangian counterparts in sufficient detail to make the parallel with
the Hamiltonian construction explicit, giving complete arguments where
the LCS overlap structure enters essentially and indicating the
relation with the Hamiltonian proof elsewhere. The logical structure is
the following: the overlap scaling of the local Lagrangians induces the
momentum scaling $T_{\alpha\beta}$ used throughout, and compatible
discrete Lagrangians generate compatible discrete right Hamiltonians,
so that the globalization theorem applies verbatim.

\paragraph{Continuous overlap structure.}

Let $\{W_\alpha\}$ be a family of open subsets of $TQ$, with
$W_\alpha$ contained in the inverse image of $V_\alpha$ under the
tangent-bundle projection, and let $L_\alpha:W_\alpha\to\R$ be a local
Lagrangian with fiber derivative
$\mathbb FL_\alpha(q,\dot q)=(q,D_2L_\alpha(q,\dot q))$. We assume that
each $W_\alpha$ is an open Legendre domain, that is, $\mathbb FL_\alpha$
is a diffeomorphism from $W_\alpha$ onto an open subset of $U_\alpha$,
and that the local inverse branches are chosen consistently on
overlaps. Regularity of $L_\alpha$ alone does not guarantee that
$W_\alpha$ can be taken to be the whole velocity fiber; this holds, for
instance, when $L_\alpha$ is hyperregular. A typical choice, when a
globally defined Lagrangian $L$ is available on compatible domains, is
$L_\alpha=e^{-\sigma_\alpha}L$. Let
$H_\alpha(q,r)=r\cdot\dot q-L_\alpha(q,\dot q)$, $r=D_2L_\alpha(q,\dot q)$,
be the associated local Hamiltonian, whose canonical Hamilton equations
are equivalent to the Euler--Lagrange equations of $L_\alpha$ on such a
domain.

\begin{proposition}
\label{prop:sec6-continuous}
Suppose that on a connected component of $W_\alpha\cap W_\beta$
\begin{equation}
    L_\beta=c_{\alpha\beta}L_\alpha.
    \label{eq:sec6-L-overlap}
\end{equation}
Then
\begin{equation}
    \mathbb FL_\beta
    =
    T_{\alpha\beta}\circ\mathbb FL_\alpha,
    \qquad
    H_\beta(q,s)
    =
    c_{\alpha\beta}H_\alpha\left(q,c_{\alpha\beta}^{-1}s\right),
    \label{eq:sec6-continuous-consequences}
\end{equation}
where the identity on the velocity variables is understood in the
first relation.
\end{proposition}

\begin{proof}
Differentiating \eqref{eq:sec6-L-overlap} in $\dot{q}$ gives
$D_2L_\beta=c_{\alpha\beta}D_2L_\alpha$, i.e., $s=c_{\alpha\beta}r$,
which is the first identity, and then
$H_\beta(q,s)=s\cdot\dot{q}-L_\beta(q,\dot{q})
=c_{\alpha\beta}[r\cdot\dot{q}-L_\alpha(q,\dot{q})]
=c_{\alpha\beta}H_\alpha(q,r)$, which is
\eqref{eq:local-Hamiltonian-overlap} after substituting
$r=c_{\alpha\beta}^{-1}s$.
\end{proof}
\noindent
Thus the multiplicative scaling of the Lagrangians is the Lagrangian
expression of the conformal cocycle \eqref{eq:conformal-cocycle}, and
it reproduces exactly the Hamiltonian overlap relation used in the
previous sections.

\paragraph{Compatible discrete Lagrangians.}

Let $L_{\alpha,h}^{E}(q_0,q_1)=\int_0^hL_\alpha(q,\dot q)\dd t$ denote
the exact discrete Lagrangian, evaluated on the solution of the
Euler--Lagrange boundary-value problem with endpoints $q_0,q_1$. Since
multiplication of a Lagrangian by a nonzero constant does not alter its
Euler--Lagrange equations, the same curve solves both boundary-value
problems on an overlap, and $L_\beta=c_{\alpha\beta}L_\alpha$ implies
$L_{\beta,h}^{E}=c_{\alpha\beta}L_{\alpha,h}^{E}$, the Lagrangian
transcription of Theorem~\ref{thm:exact-discrete-H-transformation}.

\begin{definition}
\label{def:sec6-compatible-L}
A family of discrete Lagrangians
$\{L_{\alpha,h}^{d}\}_{\alpha\in A}$ is \emph{LCS-compatible} if
\begin{equation}
    L_{\beta,h}^{d}(q_0,q_1)
    =
    c_{\alpha\beta}L_{\alpha,h}^{d}(q_0,q_1)
    \label{eq:sec6-compatible-discrete-L}
\end{equation}
on every connected component of every overlap.
\end{definition}
\noindent
As in the continuous case, multiplication of the discrete action sum
$\sum_kL^{d}_{\alpha,h}(q_k,q_{k+1})$ by the positive constant
$c_{\alpha\beta}$ leaves its stationary sequences unchanged; compatible
discrete Lagrangians therefore determine the same discrete
configuration trajectories on overlaps. Writing the left and right
local discrete Legendre transformations as
$\mathbb F^{-}L_{\alpha,h}^{d}(q_k,q_{k+1})
=(q_k,-D_1L_{\alpha,h}^{d}(q_k,q_{k+1}))$ and
$\mathbb F^{+}L_{\alpha,h}^{d}(q_k,q_{k+1})
=(q_{k+1},D_2L_{\alpha,h}^{d}(q_k,q_{k+1}))$, the associated discrete
map is
$F_{\alpha,h}^{L}=\mathbb F^+L_{\alpha,h}^{d}\circ(\mathbb F^-L_{\alpha,h}^{d})^{-1}$,
and it satisfies the following properties.

\begin{proposition}
\label{prop:sec6-discrete-consequences}
Let $\{L_{\alpha,h}^{d}\}$ be regular and LCS-compatible. Then
\begin{enumerate}[label=(\roman*)]
    \item the discrete Legendre transformations satisfy
    $\mathbb F^\pm L_{\beta,h}^{d}
    =T_{\alpha\beta}\circ\mathbb F^\pm L_{\alpha,h}^{d}$;
    \item each local map is symplectic,
    $\left(F_{\alpha,h}^{L}\right)^*\Omega_{\alpha}=\Omega_{\alpha}$;
    \item the local maps are equivariant,
    $F_{\beta,h}^{L}\circ T_{\alpha\beta}
    =T_{\alpha\beta}\circ F_{\alpha,h}^{L}$, on their common domain
    in the sense of Theorem~\ref{thm:discrete-equivariance}.
\end{enumerate}
\end{proposition}

\begin{proof}
Assertion (i) is obtained by differentiating
\eqref{eq:sec6-compatible-discrete-L} in each argument, exactly as
\eqref{eq:compatible-H-D1} was obtained from
\eqref{eq:compatible-discrete-H}; here both derivatives carry the same
factor $c_{\alpha\beta}$, so both discrete momenta scale by it.
Assertion (ii) is the transcription of
Proposition~\ref{prop:local-discrete-symplecticity}: writing
$\dd L^{d}_{\alpha,h}=\Theta^{d,+}_{\alpha,h}-\Theta^{d,-}_{\alpha,h}$
with $\Theta^{d,+}_{\alpha,h}=D_2L^{d}_{\alpha,h}(q_0,q_1)\cdot\dd q_1$
and $\Theta^{d,-}_{\alpha,h}=-D_1L^{d}_{\alpha,h}(q_0,q_1)\cdot\dd q_0$,
one more exterior derivative gives
$\dd\Theta^{d,+}_{\alpha,h}=\dd\Theta^{d,-}_{\alpha,h}$, which under
the discrete Legendre transformations is exactly
\eqref{eq:local-symplecticity-coordinate}. Assertion (iii) follows
from (i) and the definition of $F^{L}_{\alpha,h}$ by the argument of
Theorem~\ref{thm:discrete-equivariance}, with the same
connected-component proviso.
\end{proof}

\paragraph{Equivalence with the Hamiltonian construction.}

Assume that $r_1=D_2L_{\alpha,h}^{d}(q_0,q_1)$ can be solved locally
for $q_1=q_1(q_0,r_1)$ and define the associated discrete right
Hamiltonian by the discrete Legendre transform
\begin{equation}
    H_{\alpha,h}^{d+}(q_0,r_1)
    =
    r_1\cdot q_1-L_{\alpha,h}^{d}(q_0,q_1).
    \label{eq:sec6-discrete-Legendre-transform}
\end{equation}

\begin{theorem}
\label{thm:sec6-L-H}
The function \eqref{eq:sec6-discrete-Legendre-transform} satisfies
\[
    D_1H_{\alpha,h}^{d+}(q_0,r_1)
    =
    -D_1L_{\alpha,h}^{d}(q_0,q_1)
    =
    r_0,
    \qquad
    D_2H_{\alpha,h}^{d+}(q_0,r_1)
    =
    q_1,
\]
so that $L_{\alpha,h}^{d}$ and $H_{\alpha,h}^{d+}$ generate the same
phase-space map. Moreover, if the discrete Lagrangians are
LCS-compatible then so are the discrete right Hamiltonians:
\begin{equation}
    H_{\beta,h}^{d+}(q_0,s_1)
    =
    c_{\alpha\beta}
    H_{\alpha,h}^{d+}
    \left(q_0,c_{\alpha\beta}^{-1}s_1\right).
    \label{eq:sec6-H-compatibility-from-L}
\end{equation}
\end{theorem}

\begin{proof}
Differentiating \eqref{eq:sec6-discrete-Legendre-transform} and using
$r_1=D_2L^{d}_{\alpha,h}(q_0,q_1)$, the terms multiplying $\dd q_1$
cancel and
$\dd H^{d+}_{\alpha,h}=-D_1L^{d}_{\alpha,h}\cdot\dd q_0+q_1\cdot\dd r_1$,
which yields the two stated identities. For the compatibility
statement, \eqref{eq:sec6-compatible-discrete-L} gives
$s_1=D_2L^{d}_{\beta,h}(q_0,q_1)=c_{\alpha\beta}r_1$, and therefore
\[
    H^{d+}_{\beta,h}(q_0,s_1)
    =
    s_1\cdot q_1-L^{d}_{\beta,h}(q_0,q_1)
    =
    c_{\alpha\beta}
    \left[
        r_1\cdot q_1-L^{d}_{\alpha,h}(q_0,q_1)
    \right]
    =
    c_{\alpha\beta}H^{d+}_{\alpha,h}(q_0,r_1),
\]
which is \eqref{eq:sec6-H-compatibility-from-L} after substituting
$r_1=c_{\alpha\beta}^{-1}s_1$.
\end{proof}

\begin{corollary}
\label{cor:sec6-L-globalization}
A regular LCS-compatible family of discrete Lagrangians glues to a
global discrete map $F_h^{L}$ satisfying
\begin{equation}
    \left(F_h^{L}\right)^*\Omega
    =
    e^{\sigma_\alpha\circ F_h^{L}-\sigma_\alpha}\Omega.
    \label{eq:sec6-global-LCS-L}
\end{equation}
If, in addition,
$L_{\alpha,h}^{d}-L_{\alpha,h}^{E}=\cO(h^{p+1})$
in the sense of \cite[Section~2.3]{MarsWest},
with initial data in open subsets whose closures are compact and
contained in $W_\alpha$, then the resulting global
variational integrator has order $p$.
\end{corollary}

\begin{proof}
By Proposition~\ref{prop:sec6-discrete-consequences},
the local maps are symplectic and equivariant, so
Theorem~\ref{thm:global-LCS-preservation} applies.
The local order statement follows from \cite[Theorem~2.3.1]{MarsWest}.
The smooth, $h$-independent representation maps $\Psi_\alpha$
preserve the local error order, yielding the stated
order of the reconstructed map.
\end{proof}

\begin{remark}
Galerkin discrete Lagrangians provide a natural class of compatible
families: if the same approximation space and quadrature rule are used
in every local representation, then
$\mathcal L_{\beta,h}^{d}=c_{\alpha\beta}\mathcal L_{\alpha,h}^{d}$
as functions of the internal variables, the stationarity equations have
the same solutions, and
\eqref{eq:sec6-compatible-discrete-L} follows. This is the Lagrangian
counterpart of Proposition~\ref{prop:compatible-Galerkin}.
\end{remark}
\noindent
Summarizing, the scaling $L_{\beta,h}^{d}=c_{\alpha\beta}L_{\alpha,h}^{d}$
and the compatibility relation
\eqref{eq:sec6-H-compatibility-from-L} are intertwined by the discrete
Legendre transform. The Hamiltonian route remains the more general
starting point, since it requires neither the existence nor the
invertibility of a Legendre transformation; the Lagrangian route, when
available, provides an equivalent variational interpretation and an
alternative practical construction of compatible methods. Having
established both formulations, we now turn to numerical experiments
designed to test not only convergence, but also the two geometric
properties that are central to the globalization framework: overlap
compatibility and discrete conformal preservation.


\section{Numerical experiments}
\label{sec:numerics}

The experiments reported here test the three statements on which the
construction rests, separately rather than jointly:
\begin{enumerate}[label=(T\arabic*)]
    \item \label{item:T1}
    \emph{Order.} The reconstructed global map $F_h$ approximates the
    LCS flow to the order predicted by
    Proposition~\ref{prop:order-globalization} and
    Corollary~\ref{cor:GL-LCS-integrator}.
    \item \label{item:T2}
    \emph{Discrete conformal preservation.} The map $F_h$ satisfies
    $F_h^*\Omega=\rho_h\Omega$ to round-off, and not merely to the
    order of the method
    (Theorem~\ref{thm:global-LCS-preservation}).
    \item \label{item:T3}
    \emph{Overlap compatibility.} The local maps are equivariant and
    glue, so that the computed trajectory does not depend on the local
    representation in which it is computed
    (Theorem~\ref{thm:discrete-equivariance} and
    Lemma~\ref{lem:discrete-gluing}).
\end{enumerate}
\noindent
Two further experiments probe the sharpness of the framework: the
long-time behaviour of the conformal energy law
\eqref{eq:discrete-H-conformal-law}, and the failure of globalization
for families that are locally symplectic but not LCS-compatible. A
final example with a non-exact Lee form exhibits the chart switching
forced by the topology and the monodromy of the discrete conformal
factor.

\paragraph{Test problem and implementation.}

All experiments of Sections~\ref{subsec:num-convergence}
to~\ref{subsec:num-overlap} use the conformal oscillator of
Example~\ref{ex:continuous-LCS-oscillator} with $k=1$, $\mu=\tfrac12$
and initial data $(q_0,p_0,x_0,y_0)=(1,0,0,\tfrac12)$. The conformal
potential is $\sigma_\alpha=-x$, so that
$c_{\alpha\beta}=e^{C_{\alpha\beta}}$ with
$\sigma_\alpha-\sigma_\beta=C_{\alpha\beta}$ constant. In the local
chart labelled by $c=c_{\alpha\beta}$, the local Hamiltonian is, by
\eqref{eq:local-Hamiltonian-overlap},
\begin{equation}
    H_{c}(q,r,x,s)
    =
    \frac{1}{2c}e^{-x}(r^2+s^2)
    +
    c\,e^{x}\left(\frac{k}{2}q^2+\mu x\right),
    \label{eq:num-chart-Hamiltonian}
\end{equation}
which reduces to \eqref{eq:example-local-H} for $c=1$. We emphasize
that the vector field integrated in the $\beta$-representation is the
one generated by $H_\beta$, not by $H_\alpha$; this is exactly the
content of the overlap relation, and using the same local Hamiltonian
in both charts would destroy the equivariance being tested. In the LCS
variables the two-form is \eqref{eq:example-Omega-LCS-variables}, whose
matrix in the ordered basis $(\partial_q,\partial_p,\partial_x,\partial_y)$ is
\begin{equation}
    \widehat\Omega(z)
    =
    \begin{pmatrix}
        0 & 1 & p & 0\\
        -1& 0 & 0 & 0\\
        -p& 0 & 0 & 1\\
        0 & 0 & -1& 0
    \end{pmatrix}.
    \label{eq:num-Omega-matrix}
\end{equation}
Computations are carried out in double precision.
The implicit stage equations of the Gauss--Legendre methods
are solved by Newton's method with the analytic stage Jacobian,
iterating until the Newton increment satisfies
$\|\Delta Z\|_\infty\leq10\,\varepsilon_{\mathrm{mach}}\,(1+\|Z\|_\infty)$,
with at most $30$ iterations. Since this criterion is relative to the
size of the stage variables, it remains attainable in the rescaled
representations, where the momenta are of order $c_{\alpha\beta}$; an
absolute residual tolerance of order $10^{-15}$ is not attainable in
double precision for $c_{\alpha\beta}=e^{2\pi}\approx535.5$.
In exact arithmetic, Newton's method is affine covariant: with initial
guesses related by $T_{\alpha\beta}$ and compatible stopping criteria,
its iterates in the two charts correspond exactly through
$T_{\alpha\beta}$, which is consistent with the round-off-level
equivariance defects reported below.
One-step Jacobians are computed by differentiating the
implicit stage equations analytically, rather than by finite
differences. Writing $u_n=\Psi_\alpha(z_n)$ for the local
canonical coordinates and
$J_n^{\mathrm{loc}}=\partial u_{n+1}/\partial u_n$,
the Jacobian in the LCS variables is obtained by the chain rule:
$J_n
    =
    \mathrm{D}(\Psi_\alpha^{-1})(u_{n+1})
    \,J_n^{\mathrm{loc}}\,
    \mathrm{D}\Psi_\alpha(z_n).$
Here the derivative of the inverse representation map is
evaluated at the local endpoint
$u_{n+1}=\Psi_\alpha(z_{n+1})$.
This avoids finite-difference errors when measuring conformal
preservation at the level of round-off.
As a non-symplectic control, we use the classical explicit
fourth-order Runge--Kutta method applied to the same local
canonical system.

\subsection{Convergence}
\label{subsec:num-convergence}

\begin{figure}[htbp]
\centering
\includegraphics[width=0.62\textwidth]{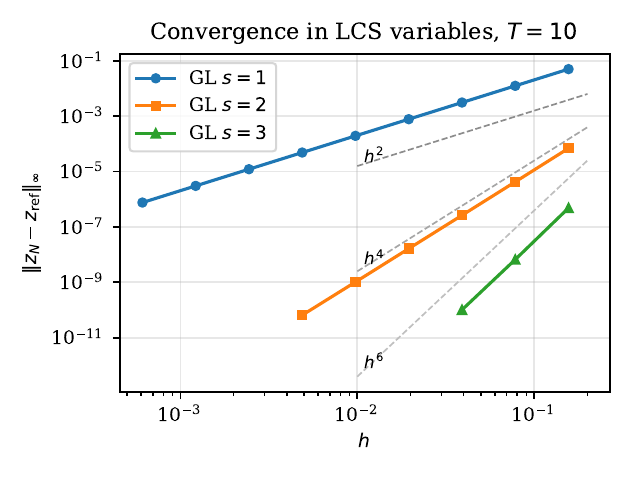}
\caption{Convergence of the reconstructed global map, confirming
orders $2$, $4$ and $6$ predicted by
Corollary~\ref{cor:GL-LCS-integrator}.}
\label{fig:convergence}
\end{figure}

We integrate to $T=10$ and measure the error
$\|z_N-z_{\mathrm{ref}}(T)\|_\infty$ in the LCS variables against a
reference trajectory computed with the three-stage Gauss--Legendre
method and $h_{\mathrm{ref}}=T/2^{19}$. Table~\ref{tab:convergence}
and Figure~\ref{fig:convergence} report the results.

\begin{table}[htbp]
\centering
\begin{tabular}{@{}lcccccc@{}}
\toprule
$h$ & $1.56\cdot10^{-1}$ & $7.81\cdot10^{-2}$ & $3.91\cdot10^{-2}$
    & $1.95\cdot10^{-2}$ & $9.77\cdot10^{-3}$ & observed order\\
\midrule
GL, $s=1$ & $4.98\cdot10^{-2}$ & $1.24\cdot10^{-2}$ & $3.10\cdot10^{-3}$
          & $7.75\cdot10^{-4}$ & $1.94\cdot10^{-4}$ & $2.00$\\
GL, $s=2$ & $6.90\cdot10^{-5}$ & $4.24\cdot10^{-6}$ & $2.63\cdot10^{-7}$
          & $1.64\cdot10^{-8}$ & $1.03\cdot10^{-9}$ & $4.01$\\
GL, $s=3$ & $4.83\cdot10^{-7}$ & $6.63\cdot10^{-9}$ & $1.01\cdot10^{-10}$
          & $^\dagger$ & $^\dagger$ & $6.11$\\
\bottomrule
\end{tabular}
\caption{Global error at $T=10$ in the LCS variables. Observed orders
are least-squares slopes over the unsaturated range. Entries marked
$^\dagger$ have reached the accuracy of the reference solution
($\approx5\cdot10^{-12}$) and are excluded from the fit.}
\label{tab:convergence}
\end{table}

\noindent
The measured orders are $2$, $4$ and $6$, in agreement with
Corollary~\ref{cor:GL-LCS-integrator}. Since the representation maps
$\Psi_\alpha$ are $h$-independent diffeomorphisms, the order measured
after reconstruction coincides with the order in the local canonical
variables, as asserted in
Proposition~\ref{prop:order-globalization}.

\subsection{Discrete conformal preservation}
\label{subsec:num-conformal}

\begin{figure}[htbp]
\centering
\includegraphics[width=0.62\textwidth]{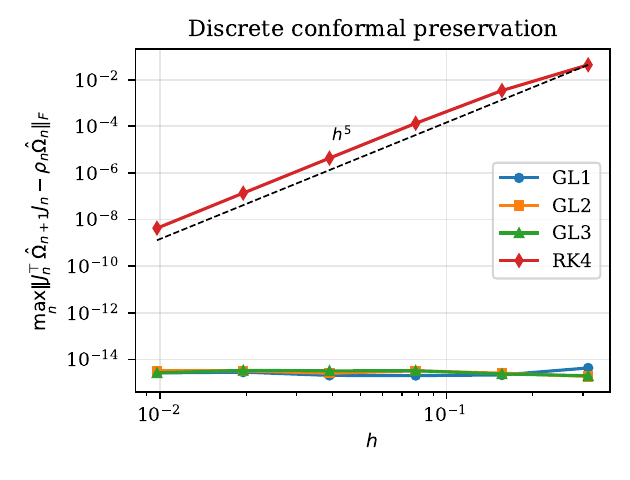}
\caption{Discrete conformal preservation. The flat curves at
$\approx3\cdot10^{-15}$ are the signature of
Theorem~\ref{thm:global-LCS-preservation}; the curve with slope
approximately five in the fine-step range is the non-symplectic
control, for which the identity holds only to the order of the
method.}
\label{fig:conformal}
\end{figure}

This is the central structural test. For each step we evaluate the
defect
\begin{equation}
    \varepsilon_n^{\mathrm{LCS}}
    =
    \left\|
        J_n^{\!\top}\widehat\Omega(z_{n+1})J_n
        -
        \rho_n\widehat\Omega(z_n)
    \right\|_F,
    \qquad
    \rho_n=e^{\sigma_\alpha(z_{n+1})-\sigma_\alpha(z_n)}
    =e^{-(x_{n+1}-x_n)},
    \label{eq:num-conformal-defect}
\end{equation}
and report $\max_n\varepsilon_n^{\mathrm{LCS}}$ over $[0,10]$ in
Table~\ref{tab:conformal} and Figure~\ref{fig:conformal}.

\begin{table}[htbp]
\centering
\begin{tabular}{@{}lcccccc@{}}
\toprule
$h$ & $3.13\cdot10^{-1}$ & $1.56\cdot10^{-1}$ & $7.81\cdot10^{-2}$
    & $3.91\cdot10^{-2}$ & $1.95\cdot10^{-2}$ & $9.77\cdot10^{-3}$\\
\midrule
GL, $s=1$ & $4.3\cdot10^{-15}$ & $2.2\cdot10^{-15}$ & $2.0\cdot10^{-15}$
          & $2.1\cdot10^{-15}$ & $2.8\cdot10^{-15}$ & $2.6\cdot10^{-15}$\\
GL, $s=2$ & $1.9\cdot10^{-15}$ & $2.6\cdot10^{-15}$ & $3.2\cdot10^{-15}$
          & $2.6\cdot10^{-15}$ & $3.2\cdot10^{-15}$ & $3.3\cdot10^{-15}$\\
GL, $s=3$ & $2.0\cdot10^{-15}$ & $2.4\cdot10^{-15}$ & $3.3\cdot10^{-15}$
          & $3.2\cdot10^{-15}$ & $3.4\cdot10^{-15}$ & $2.6\cdot10^{-15}$\\
RK4       & $4.3\cdot10^{-2}$  & $3.4\cdot10^{-3}$  & $1.3\cdot10^{-4}$
          & $4.3\cdot10^{-6}$  & $1.4\cdot10^{-7}$  & $4.3\cdot10^{-9}$\\
\bottomrule
\end{tabular}
\caption{Maximum conformal defect \eqref{eq:num-conformal-defect} over
$[0,10]$. For the Gauss--Legendre methods the defect is at round-off
and independent of $h$; for the non-symplectic RK4 control, the last
four values give an observed slope of approximately $4.96$,
consistent with an $\cO(h^5)$ one-step defect.}
\label{tab:conformal}
\end{table}

\noindent
The contrast is the point of the experiment. For the compatible
symplectic methods the defect is flat in $h$ at the level of round-off:
the identity $F_h^*\Omega=\rho_h\Omega$ holds exactly, not to the
order of the method. For the explicit Runge--Kutta control the same
quantity behaves as $\cO(h^{5})$ in the fine-step range, that is, it is
merely a consistency error. It is instructive to report also the purely
symplectic defect
$\max_n\|J_n^{\!\top}\widehat\Omega(z_{n+1})J_n-\widehat\Omega(z_n)\|_F$,
which for the implicit midpoint rule is $5.9$, $2.7$, $1.28$, $0.65$, $0.32$,
$0.16$ as $h$ is halved. This quantity is not small and decreases
like $\cO(h)$; it measures $|\rho_n-1|$ rather than any numerical
error, and confirms that the global map is conformally, not
canonically, symplectic.

\subsection{Overlap compatibility and chart switching}
\label{subsec:num-overlap}

We use two local potentials differing by a constant,
$\sigma_\alpha=-x$ and $\sigma_\beta=-x-C_{\alpha\beta}$, and the
corresponding local Hamiltonians \eqref{eq:num-chart-Hamiltonian}.
Since the Lee form of this example is exact, the two representations
cover the same domain; the experiment therefore tests the equivariance
and gluing statements directly, rather than a situation in which chart
switching is forced by the topology.

\paragraph{Equivariance defect.} For the two-stage method we evaluate
\begin{equation}
    \delta_n^{\mathrm{eq}}
    =
    \frac{\left\|
        F_{\beta,h}\left(T_{\alpha\beta}(z_n)\right)
        -
        T_{\alpha\beta}\left(F_{\alpha,h}(z_n)\right)
    \right\|_\infty}
    {\left\|T_{\alpha\beta}\left(F_{\alpha,h}(z_n)\right)\right\|_\infty}
    \label{eq:num-equivariance-defect}
\end{equation}
along the trajectory. For
$c_{\alpha\beta}\in\{e^{-1},\,2,\,e^{2\pi}\approx535.5\}$ and
$h\in\{1.56\cdot10^{-1},3.91\cdot10^{-2},9.77\cdot10^{-3}\}$, the
largest value observed over all nine runs is
$2.9\cdot10^{-16}$. For the tested transition constants, the normwise
relative equivariance defect remains at round-off level, consistently
with Theorem~\ref{thm:discrete-equivariance}.
This normalization measures the discrepancy relative to the
largest state component and does not provide a componentwise
relative error.

\paragraph{Chart switching.} We integrate with $h=T/512$ and switch
local representative at step $m\in\{1,N/4,N/2,3N/4,N-1\}$, composing
with $T_{\alpha\beta}$ at the switching step, and compare the resulting
LCS trajectory with the one computed entirely in the
$\alpha$-representation. The largest discrepancy over all switching
points is $4.1\cdot10^{-15}$ for $c_{\alpha\beta}=e^{-1}$,
$1.6\cdot10^{-15}$ for $c_{\alpha\beta}=2$ and $3.4\cdot10^{-15}$ for
$c_{\alpha\beta}=e^{2\pi}$. The glued map of
Lemma~\ref{lem:discrete-gluing} is thus independent of where the chart
is changed, to round-off.

\begin{figure}[htbp]
\centering
\includegraphics[width=0.62\textwidth]{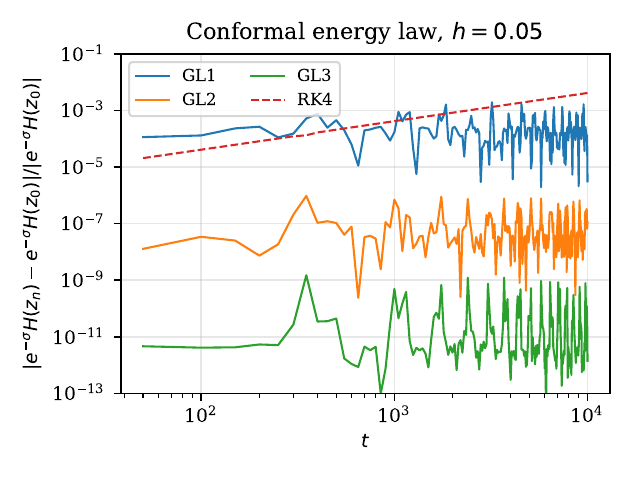}
\caption{Conformal energy law over $2\cdot10^{5}$ steps. The
Gauss--Legendre curves oscillate within a band fixed by $h$; the
non-symplectic control grows linearly.}
\label{fig:energy}
\end{figure}

\paragraph{Long-time conformal energy law.}
By Remark~\ref{rem:H-conformal-evolution}, the local Hamiltonian
$H_\alpha=e^{-\sigma_\alpha}H=e^xH$ is a first integral
of the continuous dynamics. We test the weighted energy
relation \eqref{eq:discrete-H-conformal-law} numerically
over a finite time interval by monitoring the relative error
\begin{equation}
    \Delta_n
    =
    \frac{\left|
        e^{-\sigma_\alpha(z_n)}H(z_n)
        -
        e^{-\sigma_\alpha(z_0)}H(z_0)
    \right|}
    {\left|e^{-\sigma_\alpha(z_0)}H(z_0)\right|}
    \label{eq:num-conformal-energy}
\end{equation}
over $2\cdot10^{5}$ steps with $h=0.05$, that is up to $T=10^{4}$;
Table~\ref{tab:energy} and Figure~\ref{fig:energy} show the result.

\begin{table}[htbp]
\centering
\begin{tabular}{@{}lccccc@{}}
\toprule
 & $t=T/4$ & $t=T/2$ & $t=3T/4$ & $t=T$ & $\max_n\Delta_n$\\
\midrule
GL, $s=1$ & $6.4\cdot10^{-4}$ & $2.5\cdot10^{-4}$ & $2.4\cdot10^{-5}$
          & $3.1\cdot10^{-6}$ & $2.0\cdot10^{-3}$\\
GL, $s=2$ & $5.6\cdot10^{-8}$ & $4.6\cdot10^{-8}$ & $2.4\cdot10^{-8}$
          & $1.1\cdot10^{-7}$ & $1.2\cdot10^{-6}$\\
GL, $s=3$ & $1.6\cdot10^{-11}$& $5.4\cdot10^{-11}$& $1.9\cdot10^{-12}$
          & $1.4\cdot10^{-12}$& $1.8\cdot10^{-9}$\\
RK4       & $1.0\cdot10^{-3}$ & $2.1\cdot10^{-3}$ & $3.1\cdot10^{-3}$
          & $4.1\cdot10^{-3}$ & $4.1\cdot10^{-3}$\\
\bottomrule
\end{tabular}
\caption{Relative error \eqref{eq:num-conformal-energy} in the
conformal energy law, $h=0.05$, $T=10^{4}$. The Gauss--Legendre errors
oscillate without secular growth; the explicit control drifts linearly
in $t$.}
\label{tab:energy}
\end{table}

\noindent
The qualitative distinction predicted by
Remark~\ref{rem:modified-equivariance} is clearly visible: the
Gauss--Legendre errors remain bounded and oscillatory over the whole
interval, whereas the explicit control drifts linearly, its error at
$t=T/4,T/2,3T/4,T$ being roughly proportional to $1:2:3:4$. The size of
the bounded band scales with the order of the method: repeating the
experiment to $T=10^{3}$ with $h=0.2,0.1,0.05,0.025$ gives observed
slopes $2.00$, $4.02$ and $6.04$ for $s=1,2,3$, respectively,
consistent with an $\cO(h^{2s})$ energy-error band
over the tested interval.

Across the three Gauss--Legendre runs over $[0,10^4]$,
the observed nodal values satisfy
$-3.210<x_n<0.132$ and $|q_n|<3.006$, so that the numerical orbits
remain in a bounded region over the tested interval, although the
energy level set itself is not compact.
The maximum energy errors occur near $x_n=-3.20$, where the factor
$e^{-x}$ in the kinetic term of \eqref{eq:example-local-H} is largest;
the four sampled times do not capture these peaks.

\paragraph{Locally symplectic but incompatible families.}
This experiment illustrates that local symplecticity alone does
not ensure that the local numerical maps define a single global map. We retain the compatible continuous Hamiltonians
\[
    H_\beta(q,s)
    =
    c_{\alpha\beta}
    H_\alpha(q,c_{\alpha\beta}^{-1}s),
\]
and use the two-stage Gauss--Legendre method as the underlying
discretization in both local representations. Its discrete right
Hamiltonians satisfy
\[
    H_{\beta,h}^{d+}(q,s)
    =
    c_{\alpha\beta}
    H_{\alpha,h}^{d+}(q,c_{\alpha\beta}^{-1}s).
\]
We leave $H_{\alpha,h}^{d+}$ unchanged and replace the discrete
Hamiltonian in $U_\beta$ by
\begin{equation}
    \widehat H_{\beta,h}^{d+}(q,s)
    =
    c_{\alpha\beta}
    H_{\alpha,h}^{d+}(q,c_{\alpha\beta}^{-1}s)
    +\frac{h^5}{2}\lVert s\rVert^2.
    \label{eq:incompatible-discrete-H}
\end{equation}
The power $h^5$ is chosen to match the local error order of the
fourth-order Gauss--Legendre method.

Both local maps remain symplectic. 
$\widehat H_{\beta,h}^{d+}$ defines a Type~II generating function,
and the added term does not change its mixed derivative
$D_{12}$. Thus the local nondegeneracy condition is inherited from
the unmodified discrete Hamiltonian.  The coordinate transition itself remains
$T_{\alpha\beta}(q,r)=(q,c_{\alpha\beta}r)$. To compare the two local updates, let
\[
    F_{\alpha,h}(q_k,r_k)
    =
    (q_{k+1}^{\alpha},r_{k+1}).
\]
Starting in $U_\beta$ from the corresponding point
$(q_k,s_k)$, with $s_k=c_{\alpha\beta}r_k$, the generating
relations for \eqref{eq:incompatible-discrete-H} give
\[
    s_k
    =
    c_{\alpha\beta}
    D_1H_{\alpha,h}^{d+}
    (q_k,c_{\alpha\beta}^{-1}s_{k+1}),
\]
\[
    q_{k+1}^{\beta}
    =
    D_2H_{\alpha,h}^{d+}
    (q_k,c_{\alpha\beta}^{-1}s_{k+1})
    +h^5s_{k+1}.
\]
On the locally unique solution branch, these imply
\[
    s_{k+1}=c_{\alpha\beta}r_{k+1},
    \qquad
    q_{k+1}^{\beta}
    =
    q_{k+1}^{\alpha}
    +h^5c_{\alpha\beta}r_{k+1}.
\]
Consequently, the one-step mismatch in the $\beta$ representation
is exactly
\begin{equation}
    \widehat F_{\beta,h}
    \bigl(T_{\alpha\beta}(q_k,r_k)\bigr)
    -
    T_{\alpha\beta}
    \bigl(F_{\alpha,h}(q_k,r_k)\bigr)
    =
    \bigl(h^5c_{\alpha\beta}r_{k+1},0\bigr).
    \label{eq:incompatible-map-defect}
\end{equation}
It is therefore nonzero whenever $r_{k+1}\neq0$, despite
symplecticity of both local maps. For the numerical experiment, we use the conformal oscillator
with local Hamiltonian \eqref{eq:num-chart-Hamiltonian},
taking $c=1$ in $U_\alpha$ and $c=2$ in $U_\beta$, with
$k=1$, $\mu=\tfrac12$, and initial LCS data
$(q_0,p_0,x_0,y_0)=(1,0,0,\tfrac12)$.
In the oscillator coordinates $(q,r,x,s)$, the position vector
is $(q,x)$ and the momentum vector is $(r,s)$; hence the added
term in \eqref{eq:incompatible-discrete-H} becomes
$\tfrac{h^5}{2}(r_{k+1}^2+s_{k+1}^2)$ in the $\beta$ representation.

\begin{figure}[htbp]
\centering
\includegraphics[width=0.55\textwidth]{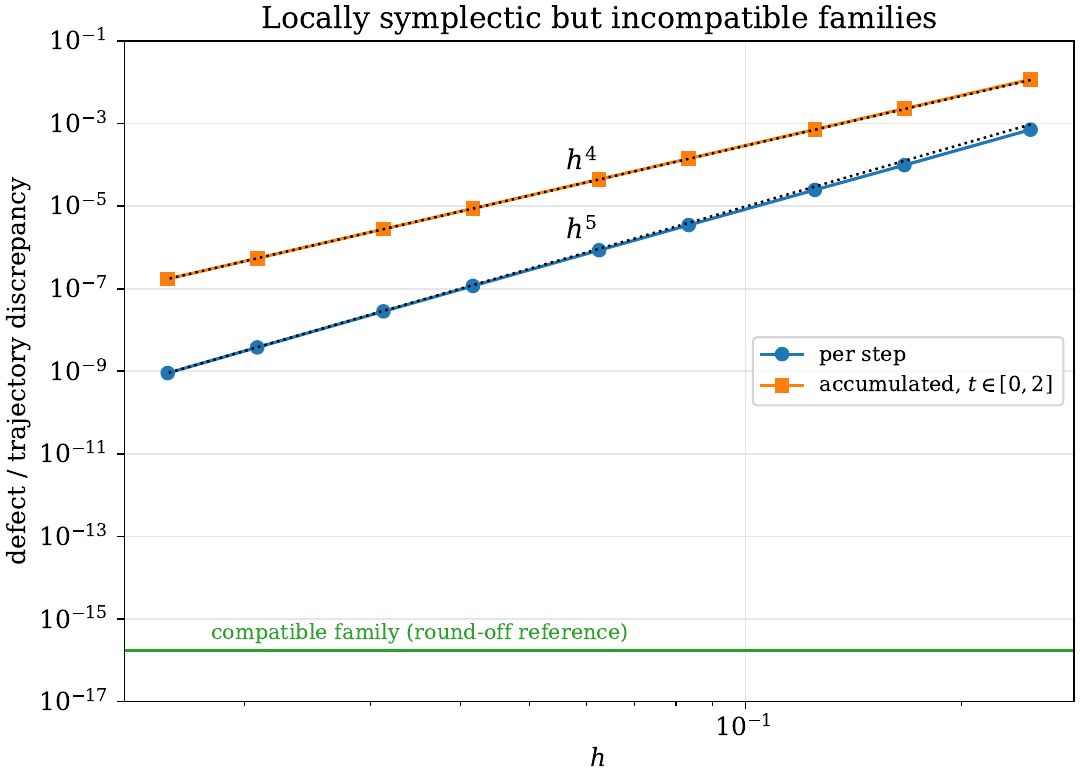}
\caption{Locally symplectic but incompatible discrete Hamiltonian
families. The equivariance defect is $\cO(h^5)$ per step, and the
trajectory discrepancy over a fixed interval is $\cO(h^4)$.
The horizontal line indicates the round-off level attained by
the compatible GL2 family.}
\label{fig:incompatible_gl2_combined}
\end{figure}

For the numerical comparison, we measure the one-step equivariance
defect defined by the norm of
\eqref{eq:incompatible-map-defect} and the maximum discrepancy
between the two independently computed trajectories, expressed
in the same LCS variables, over $t\in[0,2]$; see
Figure~\ref{fig:incompatible_gl2_combined}. The observed slopes are $4.95$ and $4.00$, respectively,
consistent with the $\cO(h^5)$ one-step defect in
\eqref{eq:incompatible-map-defect} and the $\cO(h^4)$ accumulated
discrepancy under the usual smoothness and stability assumptions.
For comparison, the unmodified local discrete maps obtained with
the Gauss--Legendre method satisfy the equivariance relation exactly.

This construction shows that a family of locally symplectic maps
can fail to define a single global numerical map when the local
updates disagree on chart overlaps. Here the disagreement is
introduced explicitly at the level of the discrete Hamiltonians,
while the continuous Hamiltonians remain compatible. Although
\eqref{eq:incompatible-map-defect} depends explicitly on
$c_{\alpha\beta}$, the added term also produces a mismatch when
$c_{\alpha\beta}=1$. Thus the experiment illustrates the failure
of overlap compatibility, rather than a defect caused solely by
a nontrivial conformal factor.

\subsection{An example with a non-exact Lee form}
\label{subsec:num-S1}

We now consider an example with a genuinely non-exact Lee form, in
order to illustrate the globalization mechanism in the presence of
nontrivial monodromy. The purpose is not to compare different
integration schemes or to repeat the convergence analysis, but to show
how locally defined discrete symplectic maps are related across a
nontrivial chart transition; we use the two-stage Gauss--Legendre
method throughout.

\begin{example}
\label{ex:S1-example}
Consider the configuration manifold $Q=S^1\times\R$ and its cotangent
bundle $M=T^*(S^1\times\R)$, with points written as $([\theta],p_\theta,x,p_x)$,
where $[\theta]\in S^1=\R/(2\pi\Z)$ denotes the angular equivalence
class. Let $V_\alpha=S^1\setminus\{[\pi]\}$ and
$V_\beta=S^1\setminus\{[0]\}$, with angular coordinates
$\theta_\alpha\in(-\pi,\pi)$ and $\theta_\beta\in(0,2\pi)$,
respectively. Using the cotangent-bundle projection
$\pi_{S^1\times\R}$, define
$U_\alpha=\pi_{S^1\times\R}^{-1}(V_\alpha\times\R)$ and
$U_\beta=\pi_{S^1\times\R}^{-1}(V_\beta\times\R)$. Let $\varphi$
denote the closed angular one-form on $S^1$, as well as its pullbacks
to $S^1\times\R$ and $M$ under the natural projections. Then
$\varphi|_{U_\alpha}=\dd\theta_\alpha$ and
$\varphi|_{U_\beta}=\dd\theta_\beta$; since $\oint_{S^1}\varphi=2\pi$,
the Lee form is closed but not globally exact, and the local Lee
potentials are $\sigma_\alpha=\theta_\alpha$ and
$\sigma_\beta=\theta_\beta$. In local angular coordinates, the
cotangent-type LCS form is represented by
\[
    \Omega
    =
    \dd\theta\wedge\dd p_\theta
    +
    \dd x\wedge\dd p_x
    +
    p_x\,\dd\theta\wedge\dd x,
\]
where $\dd\theta$ denotes the local representative of $\varphi$, and
it satisfies $\dd\Omega=\varphi\wedge\Omega$.

The overlap
$U_\alpha\cap U_\beta
    =
    \pi_{S^1\times\R}^{-1}
    \bigl((S^1\setminus\{[0],[\pi]\})\times\R\bigr)$
has two connected components. On the upper semicircle,
$0<\theta_\alpha=\theta_\beta<\pi$, so that
$\sigma_\alpha-\sigma_\beta=0$ and $c_{\alpha\beta}=1$. On the lower
semicircle, the same point of $S^1$ is represented by
$-\pi<\theta_\alpha<0$ and $\pi<\theta_\beta<2\pi$ with
$\theta_\alpha=\theta_\beta-2\pi$, so that
$\sigma_\alpha-\sigma_\beta=-2\pi$, i.e., $C_{\alpha\beta}=-2\pi$ and
$c_{\alpha\beta}=e^{-2\pi}$. Thus the two local angular coordinates
describe the same geometric point, while their local Lee potentials
differ by a nonzero constant on the second connected component of the
overlap. For $\gamma\in\{\alpha,\beta\}$ the local symplectic form is
$\Omega_\gamma=e^{-\sigma_\gamma}\Omega=\dd\theta_\gamma\wedge\dd r_\gamma
+\dd x\wedge\dd s_\gamma$ in terms of the local canonical momenta
$r_\gamma=e^{-\theta_\gamma}p_\theta$ and
$s_\gamma=e^{-\theta_\gamma}p_x$, and on the nontrivial component of
the overlap the transition between local canonical coordinates is
\[
    T_{\alpha\beta}
    (\theta_\alpha,r_\alpha,x,s_\alpha)
    =
    \bigl(
        \theta_\alpha+2\pi,
        e^{-2\pi}r_\alpha,
        x,
        e^{-2\pi}s_\alpha
    \bigr).
\]
This is an affine transition: the geometric base point is
unchanged, but its angular coordinate is translated by $2\pi$.
In the notation of \eqref{eq:affine-transition},
the base coordinates are $(\theta,x)$,
$B=I_2$, $b=(2\pi,0)^\top$, and
$c_{\alpha\beta}=e^{-2\pi}$.

We consider the globally defined Hamiltonian
$H([\theta],p_\theta,x,p_x)=p_\theta+\tfrac12(p_x^2+x^2)$, whose
definition is independent of the choice of angular lift. For
$\gamma\in\{\alpha,\beta\}$, the local Hamiltonian is
$H_\gamma=e^{-\theta_\gamma}H
=r_\gamma+\tfrac12e^{\theta_\gamma}s_\gamma^2+\tfrac12e^{-\theta_\gamma}x^2$,
with local canonical Hamilton equations
\[
    \dot\theta_\gamma=1,
    \qquad
    \dot r_\gamma
    =
    -\tfrac12e^{\theta_\gamma}s_\gamma^2
    +\tfrac12e^{-\theta_\gamma}x^2,
    \qquad
    \dot x
    =
    e^{\theta_\gamma}s_\gamma,
    \qquad
    \dot s_\gamma
    =
    -e^{-\theta_\gamma}x,
\]
which have the same form in both charts. These Hamiltonians satisfy
$H_\beta\circ T_{\alpha\beta}=e^{-2\pi}H_\alpha$.
Consequently, Proposition~\ref{prop:RK-affine-equivariance}
applies to the two-stage Gauss--Legendre method on
corresponding stage branches whose stages remain in the overlap.

For corresponding trajectories remaining in the nontrivial component
of the overlap, Proposition~\ref{prop:exact-discrete-H-affine} gives
\[
    H_{\beta,h}^{+,E}
    (\theta_{\beta,0},x_0;r_{\beta,1},s_{\beta,1})
    =
    e^{-2\pi}H_{\alpha,h}^{+,E}
    \bigl(
        \theta_{\beta,0}-2\pi,x_0;
        e^{2\pi}r_{\beta,1},e^{2\pi}s_{\beta,1}
    \bigr)
    +2\pi r_{\beta,1}.
\]
For approximate discrete right Hamiltonians, the analogous
transformation law, including the same endpoint term,
ensures compatibility by
Corollary~\ref{cor:affine-TypeII-compatibility}.

In the original LCS
variables, using a local lift of the angular coordinate, the equations
\eqref{eq:LCS-q-p-equation} read $\dot\theta=1$,
$\dot p_\theta=p_\theta-\tfrac12p_x^2+\tfrac12x^2$, $\dot x=p_x$ and
$\dot p_x=p_x-x$. For the initial condition
$(\theta_0,p_{\theta,0},x_0,p_{x,0})=(0,\tfrac12,1,0)$ the solution
of the LCS dynamics is $[\theta(t)]=[t]$,
\[
    x(t)
    =
    e^{t/2}
    \left[
        \cos\left(\frac{\sqrt3}{2}t\right)
        -
        \frac{1}{\sqrt3}
        \sin\left(\frac{\sqrt3}{2}t\right)
    \right],
\]
\[
    p_x(t)
    =
    -\frac{2}{\sqrt3}
    e^{t/2}
    \sin\left(\frac{\sqrt3}{2}t\right),
\]
and
\[
    p_\theta(t)
    =
    e^t
    -
    \frac12
    \left(
        p_x(t)^2+x(t)^2
    \right).
\]
Along this solution,
\[
    H([\theta(t)],p_\theta(t),x(t),p_x(t))
    =
    p_\theta(t)+\frac12\bigl(p_x(t)^2+x(t)^2\bigr)
    =
    e^t.
\]

\paragraph{One complete revolution and chart switching.}
Since $\dot\theta=1$, one revolution around the $S^1$ factor is
completed over $0\leq t\leq T=2\pi$. We take $N=128$ steps, so that
$h=\pi/64$, and start from $z_\alpha^0=(0,\tfrac12,1,0)$ in the local
canonical coordinates of the $\alpha$-chart (at $\theta=0$ these
coincide with the LCS variables). We switch to the $\beta$-chart at
$\theta_\alpha=\pi/2$, where $\theta_\beta=\theta_\alpha$,
$r_\beta=r_\alpha$ and $s_\beta=s_\alpha$, continue in the
$\beta$-chart until $\theta_\beta=3\pi/2$, where the nontrivial
transition $\theta_\alpha=\theta_\beta-2\pi$,
$r_\alpha=e^{2\pi}r_\beta$, $s_\alpha=e^{2\pi}s_\beta$ back to the
$\alpha$-chart is applied, and compute the final part of the
trajectory again in the $\alpha$-chart, so that $[\theta(T)]=[\theta(0)]$.
The switching locations $\pi/2$ and $3\pi/2$ are convenient choices
inside the two connected components of the overlap; the necessity of
changing local representation during a complete revolution is imposed
by the topology of $S^1$. The global momenta remain unchanged across
the nontrivial transition, since
$e^{\theta_\alpha}r_\alpha=e^{\theta_\beta}r_\beta$ and
$e^{\theta_\alpha}s_\alpha=e^{\theta_\beta}s_\beta$.

\paragraph{Monodromy of the discrete conformal factor.}
The complete revolution is precisely the situation addressed by
\eqref{eq:N-step-LCS-preservation-switching}. Since
$\theta_{n+1}=\theta_n+h$ in every chart, the chart-wise increments of
the local potentials are $\pi/2$ over the first $32$ steps in the
$\alpha$-chart, $\pi$ over the following $64$ steps in the
$\beta$-chart, and $\pi/2$ over the last $32$ steps in the
$\alpha$-chart, so that
\[
    (F_h^{128})^*\Omega
    =
    e^{2\pi}\,\Omega,
\]
although $\sigma_\alpha(z_{128})-\sigma_\alpha(z_0)=0$. This is the
discrete counterpart of
$\psi_{2\pi}^*\Omega=\exp(\oint\varphi)\,\Omega=e^{2\pi}\Omega$, and
of the conformal law $H\circ\psi_{2\pi}=e^{2\pi}H$ satisfied by the
exact solution above. Both statements are confirmed numerically. The
product $J$ of the reconstructed one-step Jacobians over the revolution
satisfies
$J^{\!\top}\widehat\Omega(z_{128})J=e^{2\pi}\widehat\Omega(z_0)$,
where $\widehat\Omega$ now denotes the matrix of $\Omega$ in the LCS
variables $(\theta,p_\theta,x,p_x)$, with a relative defect of
$1.7\cdot10^{-14}$ in the Frobenius norm, whereas the relative defect
obtained by omitting the conformal factor is
$|e^{2\pi}-1|\simeq534.492$.
Furthermore, the local energy $e^{-\theta_\gamma}H$ is conserved to a
relative accuracy of $10^{-8}$ within each chart, is continuous across
the trivial transition at $\theta=\pi/2$, and is multiplied by exactly
$e^{2\pi}$ at the nontrivial transition at $\theta=3\pi/2$, so that
$H(z_{128})/H(z_0)=e^{2\pi}\,(1+1.3\cdot10^{-8})$. The glued discrete
map therefore reproduces the monodromy of the Lee form, which no single
local potential can express.

\paragraph{One-step equivariance test.}
We next illustrate the discrete compatibility relation
$F_{\beta,h}\circ T_{\alpha\beta}=T_{\alpha\beta}\circ F_{\alpha,h}$
on the nontrivial component of the overlap. We choose a point in the
lower component of the overlap, with local canonical coordinates
$z_\alpha=(-\tfrac{\pi}{2},\tfrac12,1,\tfrac15)$
in the $\alpha$-chart, so that the transition
$z_\beta=T_{\alpha\beta}(z_\alpha)$, with $\theta_\beta=\theta_\alpha+2\pi$,
$r_\beta=e^{-2\pi}r_\alpha$ and $s_\beta=e^{-2\pi}s_\alpha$, is
nontrivial. We take $h=0.05$. We compare the state obtained by
transforming to the $\beta$-chart and then advancing one step,
$F_{\beta,h}(T_{\alpha\beta}(z_\alpha))$, with the state obtained by
advancing in the $\alpha$-chart and then transforming,
$T_{\alpha\beta}(F_{\alpha,h}(z_\alpha))$; both are expressed in the
$\beta$-representation, and we measure their discrepancy by the
normalized equivariance defect \eqref{eq:num-equivariance-defect}. For
the two-stage Gauss--Legendre scheme with $h=0.05$, we obtain
$\delta_h^{\mathrm{eq}}=3.41\times10^{-20}$. This value lies below
the unit round-off of double precision because the normalization is
dominated by $\theta_\beta\approx4.76$, whereas the momenta in the
$\beta$-representation carry the factor $e^{-2\pi}$. Componentwise,
the two routes agree exactly in $\theta$, $r_\beta$ and $x$, and differ
by approximately $2.5\times10^{-15}$ in relative terms in $s_\beta$;
the same defect measured in the $\alpha$-representation, where all
components are of order one, is $1.5\times10^{-16}$. Hence the two
discrete routes agree to round-off accuracy, consistently with the
theoretical equivariance relation. Together with the computation over
one complete revolution, this shows how the local discrete symplectic
maps are related across the nontrivial chart transition and illustrates
the globalization of the discrete dynamics in the presence of a
non-exact Lee form and nontrivial monodromy.
\end{example}


\section{Conclusions}
\label{sec:conclusions}

We have developed a Hamiltonian-first construction of
structure-preserving integrators for locally conformally symplectic
systems. The continuous LCS dynamics was first expressed through local
canonical Hamiltonian representatives and a Type~II variational
principle, without taking hyperregularity as a prerequisite for the
Hamiltonian formulation. This representation isolates the geometric
data that must be respected by a discretization: ordinary symplectic
dynamics in each local representative together with the constant
fiber-scaling cocycle on overlaps. At the discrete level, approximate
Type~II generating functions produce local symplectic maps, and the
transformation law of the exact discrete right Hamiltonian leads to a
compatibility condition for its approximations which implies
equivariance of the local discrete maps. The equivariant maps
therefore glue to a globally defined discrete dynamics on the LCS
phase space, and after reconstruction in the original variables, local
symplectic preservation becomes the discrete conformal law
\[
    F_h^*\Omega
    =
    e^{\sigma_\alpha\circ F_h-\sigma_\alpha}\Omega,
\]
which is independent of the local conformal representative, and whose
composition over several steps reproduces the periods of the Lee form
when chart switching is required.

We have exhibited two practical mechanisms for obtaining compatible
schemes: Type~II Galerkin constructions whose discrete generating
functions inherit the overlap scaling, and equivariant symplectic
Runge--Kutta methods. The Gauss--Legendre family lies in both classes,
with the implicit midpoint rule as its one-stage member. Under
regularity, the same structure is recovered from compatible discrete
Lagrangians through the discrete Legendre transformation, while the
Hamiltonian construction remains available without assuming such a
transformation globally. The numerical experiments of
Section~\ref{sec:numerics} complement these results by testing
convergence, chart switching, overlap compatibility, and the discrete
conformal preservation law, and they separate clearly the two
geometric requirements of an LCS integrator: preservation of the local
symplectic structure and compatibility of the local discretizations on
overlaps. Numerical integrators have also been developed for other
geometric and variational settings, including contact geometry; see,
for example, \cite{Simoes21, Bravetti2, Bravetti1, Vermeeren},
simple thermodynamic systems \cite{BajoLeLopez}, and non-autonomous Lagrangian systems \cite{Colombo-2023}.

The main construction uses shared base coordinates, while affine
base transitions are treated separately through the corresponding
Type~II transformation law, including the endpoint term arising
from translations. Affine equivariance also ensures compatibility
of symplectic Runge--Kutta methods under these transitions.
Extending the construction to general nonlinear base transitions
remains an open problem. The relation with conformally symplectic
integrators for damped Hamiltonian systems \cite{BhFlMo-2016, FrSuRo-2020}
also deserves emphasis: in the present LCS framework, the conformal
factor is determined by the Lee geometry rather than prescribed by a
constant damping rate. In particular, exactness of the Lee form does
not make the two settings equivalent.

\section*{Acknowledgments}

Ayten Gezici acknowledges support from The Scientific and Technological
Research Council of Türkiye (TUBITAK) under the 2219 International
Postdoctoral Research Fellowship Program, Grant No. 1059B192500320.
Melvin Leok was supported in part by the National Science Foundation
under Grants CCF-2112665, DMS-2307801, and by the Air Force Office of
Scientific Research under Grant FA9550-23-1-0279. Cristina Sard\'on
acknowledges the Programa Propio of Universidad Polit\'ecnica de Madrid
for financial support for a research stay at UCSD in the summer of 2025.

\bibliographystyle{abbrv}
\bibliography{main}

\end{document}